\documentclass[a4paper]{amsart}
\usepackage[utf8]{inputenc}
\usepackage[english]{babel}

\usepackage{amscd,amssymb,amsmath,amsthm,amstext,mathtools,esint}
\usepackage{graphicx,fancybox,epic,eepic}
\usepackage[square,numbers]{natbib}
\usepackage{booktabs}
\usepackage[hidelinks]{hyperref}
\usepackage{comment}
\usepackage{tikz,pgfplots}
\usepackage{placeins}
\usepackage{subcaption}
\usepackage{enumitem}
\usepackage{bm}

\usepackage[noabbrev, capitalise, nameinlink]{cleveref}

\usepackage{aliascnt}
\newtheorem{theorem}{Theorem}[section]

\newaliascnt{corollary}{theorem}
\newtheorem{corollary}[corollary]{Corollary}
\aliascntresetthe{corollary}

\newaliascnt{lemma}{theorem}
\newtheorem{lemma}[lemma]{Lemma}
\aliascntresetthe{lemma}

\newaliascnt{example}{theorem}

\aliascntresetthe{example}

\theoremstyle{definition}

\newaliascnt{definition}{theorem}
\newtheorem{definition}[definition]{Definition}
\aliascntresetthe{definition}

\newaliascnt{assumption}{theorem}

\aliascntresetthe{assumption}

\theoremstyle{remark}

\newaliascnt{remark}{theorem}
\newtheorem{remark}[remark]{Remark}
\aliascntresetthe{remark}

\numberwithin{theorem}{section}
\numberwithin{equation}{section}
\numberwithin{figure}{section}
\numberwithin{table}{section}

\makeatletter
\renewcommand\paragraph{\@startsection{paragraph}{4}{\z@}{.5\linespacing\@plus.7\linespacing}{-.5em}{\normalfont\bfseries}}
\makeatother

\crefname{equation}{}{}
\crefformat{equation}{\textup{#2(#1)#3}}

\AtBeginDocument{%
	\def\MR#1{}
}

\newcommand{\T}{\mathcal T}
\newcommand{\F}{\mathcal F}

\renewcommand{\forall}{\text{for all }}

\newcommand{\pw}{\mathrm{pw}}

\newcommand\Keyword[1]{{\sffamily\bfseries #1}}

\newcommand\Input{\noindent\Keyword{Input:}\enskip}

\title[High-order GLB for the Schrödinger EVP]{High-order guaranteed lower bounds for the Schrödinger eigenvalue problem}
\author[Benedikt Gr\"a\ss{}le, Moritz Hauck, and Yizhou Liang]{Benedikt Gr\"a\ss{}le$^{*}$\and Moritz Hauck$^{\dagger}$\and Yizhou Liang$^{\ddagger}$}
\thanks{$^{*}$Institut f\"ur Mathematik, Universit\"at Z\"urich, Winterthurerstrasse 190, 8057 Zurich, Switzerland}
\thanks{$^{\dagger}$Institute for Applied and Numerical Mathematics, Karlsruhe Institute of Technology, Englerstr.~2, 76131 Karlsruhe, Germany}
\thanks{$^{\ddagger}$Mathematical Institute, University of Oxford, Andrew Wiles Building, Woodstock Road, Oxford OX2 6GG, United Kingdom}
\date{}%

\begin{document}
	
	\maketitle

	\medskip 

\begin{abstract}
	We compute unconditional guaranteed lower eigenvalue bounds for Schrödinger operators using a modified hybrid high-order method with built-in correction for the potential defect.
	The method assumes exact numerical integration and allows possibly unbounded potentials in $L^p$ for $p>n/2$ and $p\geq2$ on general polyhedral domains.
	We prove that the high-order correction preserves optimal a priori convergence orders, and derive reliable and efficient a posteriori error estimators.
	Numerical experiments provide striking evidence of optimal high-order convergence rates of the adaptive algorithm.
\end{abstract}

\section{Introduction}

Guaranteed eigenvalue bounds are required whenever computed spectra enter rigorous statements, for instance in the assessment of critical loads and buckling in structural mechanics, the detection of spectral gaps, bifurcation and stability analysis, and computer-assisted existence proofs.
Combined with the upper bounds from conforming Galerkin methods for self-adjoint problems, guaranteed lower bounds provide rigorous two-sided eigenvalue enclosures.
For Schr\"odinger operators, such enclosures certify the quantum energy levels of ground and excited states.

This paper considers Schr\"odinger operators with a possibly unbounded potential $V\in L^p(\Omega)$ with $p>n/2$ and $p\geq2$ on a bounded polyhedral Lipschitz domain $\Omega\subset\mathbb{R}^n$ in $n\geq 2$ dimensions.
The corresponding Schr\"odinger eigenvalue problem seeks eigenpairs $(\lambda,u)$ with
\begin{equation} \label{eq:eigen_problem}
  \left\{
  \begin{aligned}
    -\Delta u + V u &= \lambda u && \text{in } \Omega, \\
    u &= 0 && \text{on } \partial\Omega.
  \end{aligned}
  \right.
\end{equation}
Its weak formulation seeks $(\lambda, u) \in \mathbb{R} \times H_0^1(\Omega)$ with $\|u\|_{L^2} = 1$ such that
\begin{equation}\label{eq:weak_form}
A(u, v) = \lambda b(u, v) \qquad \forall v \in H_0^1(\Omega),
\end{equation}
where $A(v, \varphi) \coloneqq (\nabla v,\nabla \varphi)_{L^2} + (Vv, \varphi)_{L^2}$ and $b(v,\varphi)\coloneqq(v,\varphi)_{L^2}$ for $v,\varphi\in H^1_0(\Omega)$.
After a possible shift of the spectrum, we may and will assume that $A(\bullet,\bullet)$ is coercive.
The normalized eigenfunctions $u$ are the critical points of the energy
\begin{align*}
	\mathcal E(v)\coloneqq A(v,v)= \int_{\Omega}\big(|\nabla v|^2+V|v|^2\big)\,\mathrm{d}x
	\qquad\text{for all }v\in H^1_0(\Omega)
\end{align*}
on the $L^2$ unit sphere in $H^1_0(\Omega)$, and each eigenvalue is the energy $\lambda=\mathcal E(u)$ of its eigenfunction.
Counting their finite multiplicity, the eigenvalues form a nondecreasing sequence $0<\lambda_1\leq\lambda_2\leq\cdots$ with $\lambda_m\to\infty$ as $m\to\infty$ 
and satisfy the Courant--Fischer min--max characterization
\begin{equation}\label{eqn:lambda_min_max}
	\lambda_m = \min_{\substack{X_m \subset H_0^1(\Omega) \\ \dim(X_m) = m}} \; \max_{v \in X_m \setminus \{0\}} \frac{\mathcal{E}(v)}{\|v\|_{L^2}^2} 
	\qquad\text{for all }m\in\mathbb{N}.
\end{equation}

Conforming Galerkin discretizations of self-adjoint elliptic eigenvalue problems provide guaranteed upper eigenvalue bounds, since minimizing over fewer subspaces in the min--max characterization~\cref{eqn:lambda_min_max} can only increase the energy.
Their convergence theory goes back to Strang and Fix~\cite{SF:AnalysisFiniteElement1973} and Babu\v{s}ka and Osborn~\cite{BabuskaOsborn1991}.
A posteriori error estimators can be used to drive adaptive mesh refinement and enable efficient higher-order conforming eigenvalue approximations; see, e.g.,~\cite{CG:RateoptimalHigherorderAdaptive2024} for a proof of rate-optimal adaptive convergence.

Guaranteed lower eigenvalue bounds (GLB) require different techniques already for the Laplace operator ($V=0$ in~\eqref{eq:eigen_problem}).
First-generation constructions based on nonconforming interpolation~\cite{Carstensen2014,Liu:FrameworkVerifiedEigenvalue2015} and mixed formulations~\cite{Gallistl2023} provide postprocessed GLB on general polyhedral domains without a priori spectral information.
These bounds are limited to low-order convergence, and their dependence on the maximal mesh size restricts the gains from local refinement.
This is particularly relevant for Schr\"odinger eigenfunctions with potential-induced localization~\cite{And58,OR:AlgorithmIdentifyingEigenvectors2023,GOP:AdaptiveRefinementEigenvalue2025}.
We refer to~\cite{Carstensen2026} for a comparison of several schemes for the Schr\"odinger problem with bounded potentials and possible simplifications for piecewise-constant potentials.

Second-generation constructions incorporate the lower-bound mechanism into the discrete problem through fine-tuned stabilization and compute direct GLB $\lambda_{h}\leq\lambda$ under a computable mesh-size condition~\cite{Carstensen2021,CP:AdaptiveGuaranteedLower2024,CGT:AdaptiveHybridHighorder2024,Tran2024}.
These schemes allow high-order variants and overcome the preceding limitations for Laplace eigenvalues: on sufficiently fine meshes, the bounds retain the approximation order of the discrete eigenvalues and recover optimal convergence rates on adaptive meshes in practice.
However, their direct extension to Schr\"odinger operators relies on orthogonality arguments that require a piecewise-constant potential, and replacing a spatially varying potential by a piecewise-constant approximation limits the accuracy of the bounds to low order.
This restriction is discussed for the extra-stabilized Crouzeix--Raviart scheme in~\cite{Carstensen2026} and also appears in the low-order bounds for the related Gross--Pitaevskii ground-state energy~\cite{GHL24,Hauck2026}.

Further approaches use explicit source-problem estimates~\cite{Liu:FrameworkVerifiedEigenvalue2015}, the Lehmann--Goerisch theorem combined with preliminary spectral bounds~\cite{LP:TwoStageFiniteElement2025}, or guaranteed a posteriori estimates~\cite{Cancs2018,CancesDussonMadayStammVohralik2020,Cha:GuaranteedStabilityBounds2026}.
Spectral Galerkin lower bounds based on explicitly known auxiliary eigensystems are developed in~\cite{Liu:GuaranteedLowerEigenvalue2026}.

\medskip\noindent\textbf{Main results.}
This paper removes the piecewise-constant restriction from second-generation schemes for the Schr\"odinger eigenvalue problem~\cref{eq:eigen_problem} with spatially varying and possibly unbounded potentials.
To this end, a hybrid high-order (HHO) discretization of the Schr\"odinger energy~$\mathcal E(\bullet)$ is modified by a negative correction that compensates the part of the potential contribution not resolved by the polynomial cell unknowns.
This correction vanishes for piecewise-constant potentials.

The focus on coercive Schrödinger forms $A(\bullet,\bullet)$ imposes no loss of generality as discussed in \Cref{rem:shift}, and simplifies the discussion of the following results.

\begin{itemize}
  \item \emph{Guaranteed lower bounds.}
  Assuming exact integration and exact algebraic eigenvalues, the discrete eigenvalues combined with a scalar postprocessing as in~\cite{CGP:ComparisonGuaranteedLower2025,Tran2024} provide GLB on arbitrary meshes without a priori spectral information.
  Once a computable condition on the maximal mesh size holds, the discrete eigenvalue $\lambda_{h,m}\leq\lambda_m$ itself is a direct GLB and local mesh refinement pays off fully.
  This condition is based on explicit control of the $L^2$ projection error.
  Potentials unbounded from below require further explicit bounds of the coercivity constant and a weighted projection error.

  \item \emph{Explicit projection bounds.}
  We derive analytical bounds for the $L^2$ projection error onto piecewise polynomials of degree $k+1$ in different (trace) norms with optimal dependence on $k$.
	  The Poincar\'e-type constant improves on the Payne--Weinberger constant already in the lowest-order case, and further leads to an explicit control for the $H^1$ stability of this projection.
  These bounds also provide the explicit parameters for the original HHO GLB~\cite{Carstensen2021} beyond the established lowest-order case $k=0$.

  \item \emph{A priori error analysis.}
  The negative correction contributes an oscillation of the potential applied to a discrete function and hence is of higher order.
  Under a mild multiplier assumption on $V$, which makes potential oscillations well defined for all of $H^1_0(\Omega)$, the modified scheme thus retains the optimal convergence rates of the underlying HHO method for eigenvalues and eigenfunctions.
  Duality arguments as in~\cite{CGT:AdaptiveHybridHighorder2024} yield improved rates for the $L^2$ and potential-weighted $L^2$ errors of the eigenfunctions.

  \item \emph{A posteriori error analysis.}
  For simple eigenvalues in $n\in\{2,3\}$ dimensions, the stabilization-free framework of~\cite{BCGT:StabilizationfreeHHOPosteriori2023,CGT:AdaptiveHybridHighorder2024} is extended to the modified scheme with spatially varying and possibly unbounded potentials.
  The residual-based estimator contains neither stabilization nor correction terms and is reliable and, up to the potential oscillation, efficient on sufficiently fine meshes.
\end{itemize}
Numerical experiments for smooth, singular, and localized eigenstates show that the adaptive algorithm driven by this a posteriori error estimator yields high-order GLB that converge with optimal rates until round-off errors dominate.

\medskip\noindent\textbf{Outline.}
\cref{sec:hho} introduces the HHO method including the potential correction term.
Explicit bounds for the $L^2$ projection and the remaining potential control precede the GLB proof in \Cref{sec:GLB}.
\Cref{sec:apriori} discusses the a priori error analysis.
The a posteriori error analysis in \cref{sec:a-posteriori} provides a reliable and efficient error estimator in $n\in \{2,3\}$ dimensions.
Numerical evidence for our theoretical results with high-order GLB on adaptive meshes concludes this paper in \Cref{sec:Numerical experiments}.

\medskip
\noindent\textbf{Notation.}
Standard notation applies for Lebesgue and Sobolev spaces and their norms and scalar products. 
For simplicity, we omit the full domain $\Omega$ in subscripts and write, for instance, $\|\bullet\|_{L^p}$, $(\bullet, \bullet)_{L^2}$, and $\|\bullet\|_{H^s}$. 
The subspace $H_0^1(\Omega)\subset H^1(\Omega)$ of functions with zero trace is equipped with the norm $\|\nabla \bullet\|_{L^2}$.

\section{Hybrid high-order method}
\label{sec:hho}

\subsection{Meshes and discrete spaces}
Consider a mesh $\mathcal{T}$ from a family of geometrically conforming and shape-regular simplicial triangulations of $\Omega$.
For each simplex element $T \in \mathcal{T}$ with the set of faces $\mathcal{F}(T)$, let $h_T \coloneqq \operatorname{diam}(T)$ denote its diameter. 
We define the piecewise constant mesh-size function $h_{\mathcal{T}}$ by $(h_{\mathcal{T}})|_T \coloneqq h_T$ for all $T \in \mathcal{T}$, and denote the maximal mesh size by $h_{\max} \coloneqq \max_{T \in \mathcal{T}} h_T$.
The set of mesh faces of $\mathcal{T}$, denoted by~$\mathcal{F}$, is partitioned into the set of interior faces $\mathcal{F}(\Omega)$ and the set of boundary faces~$\mathcal{F}(\partial \Omega)$.
The space of polynomials of total degree at most $\ell \in \mathbb{N}_0$ defined on an element $T \in \mathcal{T}$ is denoted by $\mathbb{P}^\ell(T)$. Similarly, $\mathbb{P}^k(F)$ denotes the space of polynomials of total degree at most $k \in \mathbb{N}_0$ defined on a face $F \in \mathcal{F}$. For any $T \in \mathcal{T}$, the space of broken polynomials on $\partial T$ is defined as
\[
\mathbb{P}^k(\mathcal{F}(T)) \coloneqq \bigl\{ v \in L^2(\partial T) \mid v|_F \in \mathbb{P}^k(F), \;\forall\, F \in \mathcal{F}(T) \bigr\}.
\]
The corresponding global broken polynomial spaces are
\begin{align*}
\mathbb{P}^\ell(\mathcal{T}) &\coloneqq \{v \in L^2(\Omega) \mid v|_T \in \mathbb{P}^\ell(T), \;\forall\, T \in \mathcal{T}\},\\
\mathbb{P}^k(\mathcal{F}) &\coloneqq \{v \in L^2(\Sigma) \mid v|_F \in \mathbb{P}^k(F), \;\forall\, F \in \mathcal{F}\},
\end{align*}
where $\Sigma \coloneqq \bigcup_{F \in \mathcal{F}} F$ denotes the mesh skeleton.
For all $T \in \mathcal{T}$ and all $F \in \mathcal{F}$, let $\Pi_T^\ell \colon L^2(T) \to \mathbb{P}^\ell(T)$ and $\Pi_{F}^k \colon L^2(F) \to \mathbb{P}^k(F)$ denote the $L^2$ orthogonal projections with respect to $(\bullet, \bullet)_{L^2(T)}$ and $(\bullet, \bullet)_{L^2(F)}$, respectively.
Their piecewise application defines the global $L^2$ projections
$\Pi_{\mathcal{T}}^\ell \colon L^2(\Omega) \to \mathbb{P}^\ell(\mathcal{T})$ and $\Pi_{\mathcal{F}}^k \colon L^2(\Sigma) \to \mathbb{P}^k(\mathcal{F})$.
Similarly, $\Pi_{\partial T}^k\colon L^2(\partial T)\to\mathbb P^k(\mathcal F(T))$ denotes the facewise $L^2$ projection for $T\in\mathcal{T}$.

As in the higher-order HHO constructions for GLB~\cite{Carstensen2021,Tran2024}, we use polynomial degree $k+1$ for the cell and degree $k$ for the face unknowns; the cell degree $k+1$ instead of the classical equal-order choice $k$ is essential for the GLB property, cf.~\cref{sec:GLB}.
We write $v_h=(v_{\mathcal{T}},v_{\mathcal{F}})=((v_T)_{T\in\mathcal{T}},(v_F)_{F\in\mathcal{F}})$ for $v_h\in \mathbb{P}^{k+1}(\mathcal{T}) \times \mathbb{P}^k(\mathcal{F})$
with $v_T\coloneqq v_{\mathcal{T}}|_T$, $v_F\coloneqq v_{\mathcal{F}}|_F$, and $v_{\partial T}\coloneqq(v_F)_{F\in\mathcal{F}(T)}$.
Local restrictions to $T\in\mathcal{T}$ are denoted by $v_{h,T}=(v_T,v_{\partial T})\in\mathbb{P}^{k+1}(T)\times\mathbb{P}^{k}(\mathcal{F}(T))$.
The discrete HHO space with homogeneous Dirichlet boundary conditions reads
\begin{equation*}
V_h \coloneqq \{ v_h \in \mathbb{P}^{k+1}(\mathcal{T}) \times \mathbb{P}^k(\mathcal{F}) \mid v_F = 0, \;\forall\, F \in \mathcal{F}(\partial \Omega) \}.
\end{equation*}

\subsection{Reconstruction, elliptic projection, and interpolation}

The reconstruction operator $\mathcal{R} \colon V_h \to \mathbb{P}^{k+1}(\mathcal{T})$ is defined elementwise by $(\mathcal{R} v_h)|_T \coloneqq \mathcal{R}_T v_{h,T}$ for any $v_h\in V_h$ and $T \in \mathcal{T}$, where $\mathcal{R}_T v_{h,T} \in \mathbb{P}^{k+1}(T)$ solves
\begin{align*}
(\nabla \mathcal{R}_T v_{h,T}, \nabla \varphi)_{L^2(T)} &= -(v_T, \Delta \varphi)_{L^2(T)} + (v_{\partial T}, \partial_\nu \varphi)_{L^2(\partial T)}
\quad\text{for all }\varphi \in \mathbb{P}^{k+1}(T),\\
(\mathcal{R}_T v_{h,T}, 1)_{L^2(T)} &= (v_T, 1)_{L^2(T)},
\end{align*}
with the outward normal derivative $\partial_\nu$ along $\partial T$.

Similarly, the elliptic projection operator $\mathcal{G} \colon H^1_0(\Omega) \to \mathbb{P}^{k+1}(\mathcal{T})$ is defined piecewise by $(\mathcal{G} \bullet)|_T \coloneqq \mathcal{G}_T (\bullet|_T)$, where $\mathcal{G}_T v \in \mathbb{P}^{k+1}(T)$ for any $v\in H^1(T)$ solves%
\begin{subequations}\label{eq:defGh}
    \begin{align}
        (\nabla \mathcal{G}_T v, \nabla \varphi)_{L^2(T)} &= (\nabla v, \nabla \varphi)_{L^2(T)}
        \qquad\text{for all }\varphi\in\mathbb P^{k+1}(T), \\
        (\mathcal{G}_T v, 1)_{L^2(T)} &= (v, 1)_{L^2(T)}.
    \end{align}
\end{subequations}
Denoting the piecewise gradient with respect to $\mathcal{T}$ by $\nabla_{\mathrm{pw}}$, the orthogonal projection $\mathcal{G}$ satisfies the best approximation property
\begin{equation}\label{eq:bestapproximation}
    \|\nabla_\mathrm{pw}(v - \mathcal{G} v)\|_{L^2} = \min_{v_{\mathcal{T}} \in \mathbb{P}^{k+1}(\mathcal{T})} \|\nabla_\mathrm{pw}(v - v_{\mathcal{T}})\|_{L^2}
	\qquad\text{for all }v\in H^1_0(\Omega).
\end{equation}
The interpolation operator
\begin{equation}
	\label{eq:Ih}
    \mathcal{I} \colon H^1_0(\Omega) \to V_h, \qquad \mathcal{I} v \coloneqq \big( \Pi^{k+1}_{\mathcal{T}} v,\, \Pi^k_{\mathcal{F}} v \big),
\end{equation}
connects the reconstruction and the elliptic projection through the key commuting diagram property
\begin{equation}\label{eq:keyidentity}
    \mathcal{R} \circ \mathcal{I} = \mathcal{G}.
\end{equation}

\subsection{Discrete bilinear forms}

The HHO stabilization $s_h\colon V_h\times V_h\to \mathbb{R}$ penalizes the differences between the cell and face unknowns and the corresponding values of their polynomial reconstruction.
Given a stabilization parameter $\sigma>0$, it reads
\begin{align*}
	s_h(v_h, \varphi_h) \coloneqq \sigma\sum_{T \in \mathcal{T}} s_T(v_{h,T}, \varphi_{h,T})
	\qquad\text{for all }v_h,\varphi_h\in V_h
\end{align*}
with the local contributions
\begin{align*}
s_T(v_{h,T}, \varphi_{h,T}) \coloneqq{}
&  h_T^{-1}\sum_{F \in \mathcal{F}(T)} \ell_{T,F}^{-1} \bigl( \Pi_F^k \bigl( v_{\partial T} - \mathcal{R}_T v_{h,T} \bigr),\; \Pi_F^k \bigl( \varphi_{\partial T} - \mathcal{R}_T \varphi_{h,T} \bigr) \bigr)_{L^2(F)} \\
& +  h_T^{-2} \bigl( v_T - \mathcal{R}_T v_{h,T},\; \varphi_T - \mathcal{R}_T\varphi_{h,T} \bigr)_{L^2(T)}
\qquad\text{for all }T\in\mathcal T.
\end{align*}
Following~\cite{Tran2024}, we define the face scaling parameters
\begin{equation}\label{eqn:lTF}
\ell_{T,F} \coloneqq (n+1)\frac{|F|_{n-1} \, h_T}{|T|_n}
\qquad\text{for all }F\in\mathcal F(T),
\end{equation}
where $|S|_d$ denotes the $d$-dimensional Hausdorff measure of a set $S$ for $d \in \{n-1,n\}$, abbreviated by $|S|$ if $d$ is clear.
Shape regularity implies $\ell_{T,F} \approx 1$, so that both contributions in the definition of $s_T$ scale like the squared $H^1$ seminorm.
The stabilization induces the seminorm $|v_h|_{s_h}^2 \coloneqq s_h(v_h, v_h)$ on $V_h$.

The discrete Laplace bilinear form $a_h\colon V_h\times V_h\to\mathbb R$ is defined by
\begin{equation}\label{eq:bilinear_form_ah}
    a_h(v_h, \varphi_h) \coloneqq (\nabla_{\mathrm{pw}} \mathcal{R} v_h, \nabla_{\mathrm{pw}} \mathcal{R} \varphi_h)_{L^2} + s_h(v_h, \varphi_h)
	\qquad\text{for all }v_h,\varphi_h\in V_h.
\end{equation}
It induces the discrete energy norm $\|v_h\|_{a_h}^2 \coloneqq a_h(v_h, v_h)$ on $V_h$.

Let $Q \coloneqq 1-\Pi_{\mathcal{T}}^{k+1}$ denote the $L^2$ projection error.
Given a correction parameter $\alpha>0$, the modified discrete Schr\"odinger bilinear form $A_h\colon V_h\times V_h\to\mathbb R$ reads
\begin{align}\label{eq:modified_bilinear_form}
\begin{split}
A_h(v_h, \varphi_h) &\coloneqq a_h(v_h, \varphi_h)+ (V v_{\mathcal{T}}, \varphi_{\mathcal{T}})_{L^2}
- \alpha \bigl(h_{\mathcal{T}} Q(V v_{\mathcal{T}}), h_{\mathcal{T}} Q(V \varphi_{\mathcal{T}}) \bigr)_{L^2}.
\end{split}
\end{align}
Since $V\in L^p(\Omega)\subset L^2(\Omega)$ and the cell unknowns are piecewise polynomial, $Vv_{\mathcal T}\in L^2(\Omega)$ and the correction is well defined.
If $V\in\mathbb P^0(\mathcal T)$, then
$Vv_{\mathcal T}\in\mathbb P^{k+1}(\mathcal T)$ and the correction vanishes.
The negative correction distinguishes $A_h(\bullet,\bullet)$ from standard HHO discretizations of the Schr\"odinger operator and may destroy positive definiteness on coarse meshes even for nonnegative potentials $V\geq0$.
While not relevant for the GLB result, positive definiteness always holds under the assumed coercivity of $A(\bullet,\bullet)$ on sufficiently fine meshes for fixed $\alpha>0$; see \cref{lem:shifted-coercivity} below.

The discrete mass form $b_h\colon V_h\times V_h\to \mathbb{R}$ is induced by $b(\bullet,\bullet)$ as
\begin{align*}
    b_h(v_h,\varphi_h)\coloneqq b(v_{\mathcal{T}},\varphi_{\mathcal{T}})
	\qquad\text{for all }v_h,\varphi_h\in V_h.
\end{align*}

\subsection{Discrete eigenvalue problem}
The modified HHO eigenvalue problem seeks discrete eigenpairs $(\lambda_h,u_h)\in\mathbb R\times V_h$ with $\|u_{\mathcal T}\|_{L^2}=1$ such that
\begin{equation}\label{eq:discrete_eigenproblem_modif}
A_h(u_h,\varphi_h) = \lambda_h b_h(u_h,\varphi_h) \qquad\text{for all }\varphi_h\in V_h.
\end{equation}
Since the mass, potential, and correction terms act only on the cell component, the face--face block coincides with the positive definite HHO diffusion block.
It can therefore be eliminated by a Schur complement, which reduces~\eqref{eq:discrete_eigenproblem_modif} to a symmetric eigenvalue problem for the cell unknowns as in~\cite[Sec.~3.3]{Carstensen2021}.
Hence, counting multiplicities,~\eqref{eq:discrete_eigenproblem_modif} has $N\coloneqq \dim\mathbb P^{k+1}(\mathcal{T})=|\mathcal T|\binom{n+k+1}{n}$ real discrete eigenvalues, possibly negative and enumerated in nondecreasing order
\begin{align*}
	\lambda_{h,1} \le \lambda_{h,2} \le \dots \le \lambda_{h,N}.
\end{align*}

\section{Guaranteed lower bound on eigenvalues}
\label{sec:GLB}
The modified HHO eigenvalues provide GLB via a postprocessing on arbitrary meshes and are themselves direct GLB under a computable mesh size condition.

\subsection{Explicit bounds on the $L^2$ projection}%
\label{sub:projection}
The GLB proof uses the $L^2$ orthogonal projection as HHO interpolation and requires explicit constants $C_P,C_{\mathrm{st}}>0$ in the local approximation and stability estimates
\begin{align}
	\|Qv\|_{L^2(T)} &\leq C_{P}\, h_{T} \|\nabla (v-\mathcal{G}_T v)\|_{L^2(T)},\label{eqn:CP}\\
	\|\nabla Qv\|_{L^2(T)} &\leq C_{\mathrm{st}}\|\nabla (v-\mathcal{G}_T v)\|_{L^2(T)}\label{eqn:Cst}
\end{align}
for all $T\in\mathcal T$ and $v\in H^1(T)$, where we recall $Q=1-\Pi_T^{k+1}$ on $T$. The following lemma establishes a slightly stronger form of~\eqref{eqn:CP}--\eqref{eqn:Cst}, based on the sharp bounds in~\cite{BS:ApproximationSimplicesRespect2000}. 
Define for any $k\in\mathbb{N}_0$ and with $\ell_{T,F}$ from~\eqref{eqn:lTF} the constants
\begin{align}\label{eqn:c_k}
	c_{P}^2&\coloneqq (4(k+2)(k+n+2))^{-1} \\
	c_{\mathrm{tr}}^2&\coloneqq c_P\bigl(c_P+2/(n+1)\bigr) \\
	c_{\mathrm{st}}^2&\coloneqq 1+
		(k+1)(k+n)(n(n+1))^{-1}
		c_{\mathrm{tr}}^2\sum\nolimits_{F\in\mathcal{F}(T)}\ell_{T,F}^2\label{eqn:c_st}.
\end{align}%
Let $\Pi_T^k$ act componentwise on gradients.
\begin{lemma}[Explicit projection bounds]\label{lem:explicit_projection_bound}
	Any nondegenerate simplex $T\subset\mathbb{R}^n$ and any $v\in H^1(T)$ satisfy
	\begin{align}
		\|Qv\|_{L^2(T)}^2
		&\leq c_{P}^2h_T^2\|(1-\Pi_T^k)\nabla v\|_{L^2(T)}^2,\label{eqn:proj_L2}\\
		\sum\nolimits_{F\in\mathcal{F}(T)} \ell_{T,F}^{-1}\|Qv\|_{L^2(F)}^2
		&\leq c_{\mathrm{tr}}^2h_T\|(1-\Pi_T^k)\nabla v\|_{L^2(T)}^2,\label{eqn:proj_trace}\\
		\|\nabla Qv\|_{L^2(T)}^2
		&\leq c_{\mathrm{st}}^2\|(1-\Pi_T^k)\nabla v\|_{L^2(T)}^2.\label{eqn:proj_H1}
	\end{align}
\end{lemma}

Since $\Pi_T^k$ is the best approximation in $L^2(T)$, the estimates~\eqref{eqn:proj_L2} and~\eqref{eqn:proj_H1} imply~\eqref{eqn:CP}--\eqref{eqn:Cst} with $C_P\leq c_P$ and $C_{\mathrm{st}}\leq c_{\mathrm{st}}$.
The constant $c_{P}$ scales optimally like $(k+1)^{-1}$ and improves on the Payne--Weinberger constant $1/\pi$ for all $k\geq0$.
\Cref{lem:explicit_projection_bound} also provides explicit expressions for the constants in~\cite[Thm.~2.3]{Carstensen2021} of the original HHO GLB, which were only available in the lowest-order case $k=0$ before; the growth $\sqrt{k+1}$ of $c_{\mathrm{st}}$ is optimal~\cite[Lem.~2.1]{CGT:AdaptiveHybridHighorder2024} for $n\in\{2,3\}$.%

\begin{proof}[Proof of \Cref{lem:explicit_projection_bound}]
	Let $P_1,\dots,P_{n+1}$ denote the vertices of $T$ with barycentre $x_T$, outer unit normal $\nu_T$, and barycentric coordinates $\lambda_1,\dots,\lambda_{n+1}$.
	Set
	\begin{align*}
		e\coloneqq Qv,
		\qquad
		q\coloneqq \Pi_{T}^k\nabla v-\nabla \Pi_T^{k+1}v
		\quad\text{so that}\quad
		r\coloneqq\nabla e-q=(1-\Pi_{T}^k)\nabla v.
	\end{align*}

	\medskip\noindent\emph{Proof of~\eqref{eqn:proj_L2}}.
	This is a consequence of~\cite[Thm.~3.3]{BS:ApproximationSimplicesRespect2000}, which provides
	\begin{align}\label{eqn:BS_bound}
		\|e\|_{L^2(T)}=\|Qe\|_{L^2(T)}\leq 2c_P\|M^{1/2}\nabla e\|_{L^2(T)}
	\end{align}
	with the symmetric matrix-valued function $M\in \mathbb P^{2}(T)^{n\times n}$ given by
	\begin{align*}
		M(x) = \sum^{}_{1\leq i<j\leq n+1} \lambda_{i}(x)\lambda_{j}(x)(P_i-P_j) (P_i-P_j)^\top
		\qquad\text{for all }x\in T.
	\end{align*}
	Since $s_i\coloneqq\xi\cdot P_i-\min_{1\leq j\leq n+1}\xi\cdot P_j\leq h_T$ holds for any $\xi\in \mathbb{R}^n$ with $|\xi|=1$ and $i=1,\dots,n+1$, a straightforward calculation with $s=\sum_{i=1}^{n+1}\lambda_is_i$ shows
	\begin{align*}
		\xi\cdot M \xi= \hspace{-.5em}\sum^{}_{1\leq i<j\leq n+1} \hspace{-.5em}\lambda_i\lambda_j(s_i-s_j)^2
				= \sum^{n+1}_{i=1} \lambda_i s_i^2 - s^2
				\leq s(h_T-s)= \frac{h_T^2}{4}-\left(\frac{h_T}{2}-s\right)^2.
	\end{align*}
	Hence the spectral norm of $M^{1/2}$ is bounded by $h_T/2$.
	Since $M\nu_{T}\equiv 0$ on $\partial T$, an integration by parts and the orthogonality of $e$ to $\operatorname{div}(Mq)\in \mathbb{P}^{k+1}(T)$ provide
	\begin{align*}
		(M\nabla e, q)_{L^2(T)}=-(e, \operatorname{div}(Mq) )_{L^2(T)}=0.
	\end{align*}
	This orthogonality and the spectral norm bound reveal
	\begin{align*}
		\|M^{1/2}\nabla e\|_{L^2(T)}^2 = (M\nabla e, r)_{L^2(T)}
		\leq (h_T/2) \|M^{1/2}\nabla e\|_{L^2(T)}\|r\|_{L^2(T)}.
	\end{align*}
	The combination with~\eqref{eqn:BS_bound} concludes the proof of~\eqref{eqn:proj_L2}.

	\medskip
	\noindent\emph{Proof of~\eqref{eqn:proj_trace}}.
	Following the proof of~\cite[Lem.~3]{BCGT:StabilizationfreeHHOPosteriori2023},
	the trace identity on the sub-simplices joining the faces of $T$ with its barycentre $x_{T}$ results in
	\begin{align*}
		\sum_{F\in\mathcal F(T)}
		\frac{n|T|}{(n+1)|F|}\|e\|_{L^2(F)}^2  
		&= n \|e\|_{L^2(T)}^2 + 2 (e, (x-x_{T})\cdot \nabla e)_{L^2(T)}\\
		&= n \|e\|_{L^2(T)}^2 + 2 (e, (x-x_{T})\cdot r)_{L^2(T)}
	\end{align*}
	with the orthogonality of $e$ to $(x-x_{T})\cdot q\in \mathbb{P}^{k+1}(T)$ in the last step.
	Hence the geometric bound $\max_{x\in T}|x-x_{T}|<nh_T/(n+1)$ reveals the multiplicative form
	\begin{align*}
	\sum_{F\in\mathcal F(T)}
		\ell_{T,F}^{-1}h_T\|e\|_{L^2(F)}^2
		\leq\|e\|_{L^2(T)}^2
		+\frac{2}{n+1}h_T\|e\|_{L^2(T)}\|r\|_{L^2(T)}.
	\end{align*}
	This and~\eqref{eqn:proj_L2} conclude the proof of~\eqref{eqn:proj_trace}.

	\medskip
	\noindent\emph{Proof of~\eqref{eqn:proj_H1}}.
	Orthogonality and an integration by parts verify
	\begin{align*}
		\|q\|_{L^2(T)}^2 = (q, \nabla e)_{L^2(T)} 
		= (q\cdot \nu_T, e)_{L^2(\partial T)}.
	\end{align*}
	The inverse estimate~\cite[Thm.~5]{WH:ConstantsHpfiniteElement2003} applied componentwise on $F\in\mathcal{F}(T)$ shows
	\begin{align*}
		\frac{n}{(k+1)(k+n)}
		\|q\|_{L^2(F)}^2
		\leq \frac{|F|}{|T|}
		\|q\|_{L^2(T)}^2
		\leq \frac{\ell_{T,F}}{(n+1)h_T}
		\|q\|_{L^2(T)}^2.
	\end{align*}
	Together with~\eqref{eqn:proj_trace}, this and Cauchy inequalities imply
	\begin{align*}
		\|q\|_{L^2(T)}
		\leq c_{\mathrm{tr}}
		\left(
		\frac{(k+1)(k+n)}{n(n+1)}
		\sum_{F\in\mathcal{F}(T)}\ell_{T,F}^2
		\right)^{1/2}
		\|r\|_{L^2(T)}.
	\end{align*}
	This and the orthogonal split $\|\nabla e\|_{L^2(T)}^2=\|r\|_{L^2(T)}^2 + \|q\|_{L^2(T)}^2$ yield~\eqref{eqn:proj_H1}.
\end{proof}
It is interesting to note that \Cref{lem:explicit_projection_bound} also implies an explicit bound for the $H^1$ stability constant of the local $L^2$ projection.
\begin{corollary}[Explicit $H^1$ stability of $\Pi_{T}^{k+1}$]\label{cor:H1_stability}
	Any simplex $T\subset\mathbb{R}^n$ and $v\in H^1(T)$ satisfy $\|\nabla\Pi_T^{k+1}v\|_{L^2(T)} \leq c_{\mathrm{st}}\|\nabla v\|_{L^2(T)}$ with $c_{\mathrm{st}}$ from~\eqref{eqn:c_st}.
\end{corollary}
\begin{proof}
	This follows as in~\cite[Lem.~2.1]{CGT:AdaptiveHybridHighorder2024} from Kato's lemma, i.e., equal operator norms of $\Pi_{T}^{k+1}$ and $Q$ in the Hilbert space $H^1(T)/\mathbb R$ with norm $\|\nabla\bullet\|_{L^2(T)}$.
\end{proof}

\subsection{Explicit constants for the modified scheme}%
\label{sub:Certification constants}
The GLB proof requires explicit information in the form of a stabilization bound $\Lambda>0$ and, if $V$ takes negative values, a coercivity bound $c_A>0$ and a potential-weighted projection bound $\delta_h\geq0$.

\paragraph{Stabilization}
The following lemma provides the stabilization bound $\Lambda>0$.
\begin{lemma}[Stabilization estimate]\label{lem:stability}
	Any $T\in\mathcal T$ and $v\in H^1(T)$ with local interpolation $v_{h,T}\coloneqq(\Pi_T^{k+1}v,\Pi_{\partial T}^k v)$ satisfy
\begin{align}\label{eq:proof_stabilization}
	s_T(v_{h,T},v_{h,T}) \leq \Lambda \|\nabla(v - \mathcal{G}_T v)\|_{L^2(T)}^2
	\qquad\text{with }\Lambda\leq 2/\pi^2 + 2/(n\pi).
\end{align}
\end{lemma}
\begin{proof}
	This is known from~\cite{Tran2024} as a consequence of Lemma 3.1 therein.
\end{proof}
The optimal local values of $\Lambda$, $C_P$, and $C_{\mathrm{st}}$ depend only on $k$ and the shape of $T$.
They are characterized by local variational eigenvalue problems and can be enclosed using the verified eigenvalue framework of~\cite{Liu:FrameworkVerifiedEigenvalue2015}, with its extension to positive semidefinite bilinear forms in~\cite[Thm.~2.4]{YXL:GuaranteedEigenvalueBounds2019}.
\Cref{tab:Lambda} below lists numerical approximations for right-isosceles triangles.

\paragraph{Coercivity bound}
The postprocessing involves $c_A>0$ in the coercivity estimate
\begin{align}\label{eqn:Garding}
	A(v,v) \geq c_A \|\nabla v\|_{L^2}^2
	\qquad\text{for all }v\in H^1_0(\Omega).
\end{align}
For nonnegative $V\geq0$, the natural choice is $c_A=1$.
For general $V\in L^p(\Omega)$, split $V=V_+-V_-$ into $V_{\pm}=(\pm V)_+$ with $(\bullet)_+=\max\{0,\bullet\}$.
It suffices to consider $V_-$ with small $L^p$ norm after a possible spectral shift as discussed in \Cref{rem:shift} below.
The Sobolev embedding then absorbs this perturbation.

Explicit analytical bounds on the Sobolev constant $C_{S,r}>0$ for $r=2p/(p-1)$~in 
\begin{equation}\label{eq:global_Dirichlet_Sobolev}
	\|v\|_{L^r}\leq C_{S,r}\|\nabla v\|_{L^2}
    \qquad\text{for all }v\in H_0^1(\Omega)
\end{equation}
follow as in~\cite[Cor.~A.2]{TSMO:SharpNumericalInclusion2017} by zero extension and H\"older's inequality from the critical embedding $W^{1,q}(\mathbb R^n)\hookrightarrow L^r(\mathbb R^n)$ for ${q=nr/(n+r)}$ with optimal constants from~\cite{Aub:ProblemesIsoperimetriquesEspaces1976,Talenti1976}.
Alternative bounds follow from an interpolation inequality whose $L^2$ factor is controlled by a lower bound for the first Dirichlet Laplace eigenvalue~\cite[Lem.~7.10]{NPW:NumericalVerificationMethods2019}.
On bounded convex domains, the numerical algorithm in~\cite[Thm.~1.1]{TSMO:SharpNumericalInclusion2017} also provides sharp verified enclosures of $C_{S,r}$.
\begin{lemma}[Explicit coercivity bound]\label{lem:Garding_Lp}
	If $C_->0$ bounds $\|V_-\|_{L^p}\leq C_-$ and satisfies $C_{S,r}^2 C_-<1$, then~\eqref{eqn:Garding} holds with $c_A=1-C_{S,r}^2 C_->0$.
\end{lemma}
\begin{proof}
	This follows from a Hölder inequality for $A(v,v)\geq\|\nabla v\|_{L^2}^2 - (V_- v,v)_{L^2}$ with any $v\in H^1_0(\Omega)$ and $(V_-v,v)_{L^2}\leq C_{S,r}^2\|V_-\|_{L^p}\|\nabla v\|_{L^2}^2$ using~\eqref{eq:global_Dirichlet_Sobolev}.
\end{proof}

\paragraph{Potential-weighted projection error}
The direct GLB condition further involves the constant $\delta_h\geq0$ in the potential-weighted estimate
\begin{align}\label{eqn:V_negative_control}
	\|V_-^{1/2}Qv\|_{L^2}
	\leq \delta_h\|\nabla_\pw(v-\mathcal Gv)\|_{L^2}
	\qquad\text{for all }v\in H^1_0(\Omega).
\end{align}
For nonnegative $V\geq0$, the natural choice is $\delta_h=0$.
Unlike $C_P$, the optimal $\delta_h$ depends on $V_-$ and \emph{cannot} be tabulated from element shapes if $V_-$ varies in space.
For $V_-\in L^\infty(\Omega)$ with $\|V_-\|_{L^\infty}\leq C_-$, one may take $\delta_h=C_Ph_{\max}\sqrt{C_-}$.

For $V_-\in L^p(\Omega)$ with $p<\infty$, an explicit $\delta_h$ follows from a constant $C_{P,r}>0$ with $r=2p/(p-1)$ such that all $T\in\mathcal T$ and $v\in H^1(T)$ satisfy
\begin{align}\label{eqn:Poincare_s}
\|Qv\|_{L^r(T)} \leq C_{P,r}h_T^{\,1+n/r-n/2}
\|\nabla(v-\mathcal G_Tv)\|_{L^2(T)}.
\end{align}
\begin{lemma}[Explicit negative-potential bound]\label{lem:negative_potential_Lp}
	If $C_->0$ bounds $\|V_-\|_{L^p}\leq C_-$, then~\eqref{eqn:V_negative_control} holds with $\delta_h^2\coloneqq C_{P,r}^2 h_{\max}^{2-n/p}C_-$.
\end{lemma}
\begin{proof}
	Hölder's inequality and~\eqref{eqn:Poincare_s} show, using $2-n/p>0$, for any $T\in\mathcal T$ that
	\begin{align*}
		\|V_-^{1/2}Qv\|_{L^2(T)}^2
		\leq\|V_-\|_{L^p(T)}\|Qv\|_{L^r(T)}^2
		\leq C_-C_{P,r}^2h_{\max}^{2-n/p}\|\nabla(v-\mathcal G_Tv)\|_{L^2(T)}^2
	\end{align*}
	The sum over $T\in\mathcal T$ concludes the proof.
\end{proof}

Like $C_P=C_{P,2}$, the optimal constant $C_{P,r}$ depends only on $k$ and the shape of~$T$, but its verified computation is significantly more involved for $r>2$, since the associated extremal problem is nonlinear and already convergence may be an issue.

Analytical bounds may combine~\eqref{eqn:Cst} with constants $D_s>0$ for $s\in I_r$ with $I_r\coloneqq[r,\infty)$ for $n=2$ and $I_r\coloneqq [r,2n/(n-2)]$ for $n\geq3$ in
\begin{align}\label{eqn:Ds}
    \|v-\Pi_T^0v\|_{L^s(T)}
	&\leq D_s h_T^{1+n/s-n/2}\|\nabla v\|_{L^2(T)}
    \qquad\text{for all }v\in H^1(T).
\end{align}
Explicit analytical bounds~\cite[Thm.~3.2--3.3]{MTSO:EstimationSobolevEmbedding2017} for $D_s$ on bounded convex domains are available.
Recall $C_P\leq c_P$ and $C_{\mathrm{st}}\leq c_{\mathrm{st}}$ from~\Cref{lem:explicit_projection_bound}.
\begin{lemma}[$L^r$ projection error]\label{lem:projection_sobolev}
	The estimate~\eqref{eqn:Poincare_s} holds for
	\begin{align*}
		C_{P,r}
		\coloneqq
		\inf_{s\in I_r}
		C_P^{\theta_s}(C_{\mathrm{st}}D_s)^{1-\theta_s},
		\qquad\text{where }
		\theta_s\coloneqq
		\frac{\tfrac1r-\tfrac1s}{\tfrac12-\tfrac1s}.
\end{align*}
\end{lemma}
\begin{proof}
    Let $s\in I_r$ and $\theta=\theta_s$. Hölder's inequality for $|Qv|=|Qv|^{\theta}|Qv|^{1-\theta}$ with $1/r=\theta/2+(1-\theta)/s$, \eqref{eqn:CP}, and~\eqref{eqn:Ds} for $Qv=(1-\Pi_T^0)Qv$ with~\eqref{eqn:Cst} show
    \begin{align*}
    	\|Qv\|_{L^r(T)}
    	&\leq\|Qv\|_{L^2(T)}^{\theta}\|Qv\|_{L^s(T)}^{1-\theta}\\
    	&\leq(C_Ph_T)^{\theta}\bigl(C_{\mathrm{st}}D_sh_T^{1+n/s-n/2}\bigr)^{1-\theta}\|\nabla(v-\mathcal G_Tv)\|_{L^2(T)}\\
    	&=C_P^{\theta}(C_{\mathrm{st}}D_s)^{1-\theta}h_T^{1+n/r-n/2}\|\nabla(v-\mathcal G_Tv)\|_{L^2(T)}
    \end{align*}
    with $\theta+(1-\theta)(1+n/s-n/2)=1+n/r-n/2$ in the last step. The infimum over $s\in I_r$ concludes the proof.
\end{proof}

Alternative bounds for $C_{P,r}$ follow from explicit Gagliardo--Nirenberg inequalities, which interpolate the $L^r$ norm between the $L^2$ norm and the $H^1$ seminorm.
On right-isosceles triangles, reflection reduces this to the periodic setting, where the inequalities in~\cite{IKZ:ApplicationsLiebThirring2022,ZI:ClassInterpolationInequalities2023} provide the following simple bound.
\begin{lemma}[$L^r$ projection error for the right-isosceles triangle]\label{lem:Lr_right_iso}
	Any right-isosceles triangle $T\subset\mathbb{R}^2$ satisfies~\eqref{eqn:Poincare_s} for $2\leq r<\infty$ with
	\begin{align*}
		C_{P,r}=\sqrt{\frac{r}{2}} \left(\frac{2}{\pi}\right)^{1/2-1/r}C_{P}^{2/r}C_{\mathrm{st}}^{1-2/r}.
	\end{align*}
\end{lemma}
\begin{proof}
	Let $v\in H^1(T)$. 
	After a rigid motion, $T=\operatorname{conv}\{0,ae_1,ae_2\}$ with $a>0$.
	The even reflection across the hypotenuse extends $Qv$ to $(0,a)^2$, and the even extension in both coordinates to $(-a,a)^2$, which consists of eight copies of~$T$, is $2a$-periodic and defines $w\in H^1(\mathbb{T}^2)$ on $\mathbb{T}^2=\mathbb{R}^2/(2a\mathbb{Z})^2$.
	Hence $\|w\|_{L^s(\mathbb{T}^2)}^s=8\|Qv\|_{L^s(T)}^s$ for all $s<\infty$ and $\|\nabla w\|_{L^2(\mathbb{T}^2)}^2=8\|\nabla Qv\|_{L^2(T)}^2$.
	As $w$ has zero mean,~\cite[(1.8)]{ZI:ClassInterpolationInequalities2023} with constant $\sqrt{r/2}\,(4\pi)^{1/r-1/2}$ and the factor $8^{1/2-1/r}$ from the eight copies~provides
	\begin{align*}
		\|Qv\|_{L^r(T)}
		\le
		\sqrt{\frac r2}
		\left(\frac2\pi\right)^{1/2-1/r}
		\|Qv\|_{L^2(T)}^{2/r}
		\|\nabla Qv\|_{L^2(T)}^{1-2/r}.
	\end{align*}
	This and~\eqref{eqn:CP}--\eqref{eqn:Cst} conclude the proof.
\end{proof}%
\subsection{Guaranteed lower bounds}%
\label{sub:Guaranteed lower bounds}
The parameters $\sigma$ and $\alpha$ enter with the constants $\Lambda$, $C_P$, and $\delta_h$ from \Cref{sub:projection,sub:Certification constants} into the certification constant
\begin{align}\label{eqn:gamma_def}
	\gamma(\lambda_{h,m})\coloneqq (\sigma\Lambda + C_P^2(\alpha^{-1} + h_{\max}^2(\lambda_{h,m})_+) + \delta_h^2-1)_{+}.
\end{align}
The coercivity bound $c_A>0$ only contributes to the postprocessing if $\gamma(\lambda_{h,m})>0$.
\begin{theorem}[Guaranteed lower bound]\label{theo:main_lower_bound}
	For any mesh $\mathcal T$, parameters $\sigma,\alpha>0$, and $1\leq m\leq N$, the $m$-th eigenvalues $\lambda_m$ of~\eqref{eq:weak_form} and $\lambda_{h,m}$ of~\cref{eq:discrete_eigenproblem_modif} satisfy
\begin{equation}
\mathrm{GLB}(m) \coloneqq \frac{c_A(\lambda_{h,m})_+}{c_A+\gamma(\lambda_{h,m})} \leq \lambda_m.
\end{equation}
In particular, if
\[
\sigma\Lambda + C_P^2(\alpha^{-1} + h_{\max}^2(\lambda_{h,m})_+) + \delta_h^2 \leq 1,
\]
then $\gamma(\lambda_{h,m}) = 0$ and $\lambda_{h,m} \leq \operatorname{GLB}(m)\leq \lambda_m$ is a direct lower eigenvalue bound.
\end{theorem}

\begin{proof}
	Let $\phi_1,\dots,\phi_m$ denote $L^2$-orthonormal eigenfunctions of~\eqref{eq:weak_form} associated with eigenvalues $\lambda_1,\dots,\lambda_m$ that define $\mathbb{E}\coloneqq \operatorname{span}\{\phi_1,\dots,\phi_m\}$. 
	The Schrödinger energy $\mathcal{E}(\phi)=A(\phi,\phi)$ satisfies $\mathcal{E}(\phi)\leq \lambda_m$ for any $\phi\in\mathbb E$ with $\|\phi\|_{L^2}=1$ by the $A$- and $L^2$-orthogonality of the eigenfunctions.
	The best approximation~\eqref{eq:bestapproximation} of $\mathcal G$ and coercivity~\eqref{eqn:Garding}~imply
	\begin{equation}\label{eqn:bestapprox_bound}
	\|\nabla_\pw(\phi-\mathcal G\phi)\|_{L^2}^2
	\leq \|\nabla\phi\|_{L^2}^2
	\leq c_A^{-1}\mathcal E(\phi)
	\leq c_A^{-1}\lambda_m.
	\end{equation}
	This proof distinguishes between linearly dependent and independent cell projections $\Pi_{\mathcal{T}}^{k+1}\phi_1,\dots,\Pi_{\mathcal{T}}^{k+1}\phi_m$.

\medskip\noindent\emph{Case 1: linearly dependent cell projections.} Then some $\phi\in \mathbb{E}$ with $\|\phi\|_{L^2}=1$ satisfies $\Pi_{\mathcal{T}}^{k+1}\phi=0$, and the Poincaré inequality~\eqref{eqn:CP} and~\cref{eq:bestapproximation} yield
\begin{align*}
	1 = \|\phi\|_{L^2}^2 = \|Q\phi\|_{L^2}^2 \leq h_{\max}^2 C_P^2\|\nabla\phi\|_{L^2}^2.
\end{align*}
The definition~\eqref{eqn:gamma_def} of $\gamma(\lambda_{h,m})$ provides
$h_{\max}^2C_P^2(\lambda_{h,m})_+\leq1-\delta_h^2+\gamma(\lambda_{h,m})$.
Since $Q\phi=\phi$, the potential-weighted estimate~\eqref{eqn:V_negative_control} and~\eqref{eq:bestapproximation} imply
$\|V_-^{1/2}\phi\|_{L^2}^2\leq\delta_h^2\|\nabla\phi\|_{L^2}^2$.
Hence $(1-\delta_h^2)\|\nabla\phi\|_{L^2}^2\leq\mathcal E(\phi)$, so that
\begin{align*}
	\lambda_{h,m}
	&\leq(\lambda_{h,m})_+
	\leq\bigl(1-\delta_h^2+\gamma(\lambda_{h,m})\bigr)
	\|\nabla\phi\|_{L^2}^2
	\leq\mathcal E(\phi)
	+\gamma(\lambda_{h,m})\|\nabla\phi\|_{L^2}^2.
\end{align*}
This, $\mathcal{E}(\phi)\leq \lambda_m$, and $\|\nabla \phi\|_{L^2}^2\leq \lambda_m /c_A$ from~\eqref{eqn:bestapprox_bound} conclude this case.

\medskip\noindent\emph{Case 2: linearly independent cell projections.} 
Here
\begin{align*}
	\mathbb{E}_h\coloneqq \operatorname{span}\{\mathcal I\phi_1,\dots,\mathcal I\phi_m\}\subset V_h
\end{align*}
is an $m$-dimensional subspace with $v_{\mathcal T}\neq0$ for all $v_h\in\mathbb E_h\setminus\{0\}$. Hence the discrete min--max principle, also known as the Courant--Fischer theorem, reveals
\begin{align*}
	\lambda_{h,m}= \min_{\substack{X_{h} \subset V_h,\ \dim(X_{h})=m\\ v_{\mathcal T}\neq0\text{ for all }v_h\in X_h\setminus\{0\}}} \; \max_{v_h \in X_{h} \setminus \{0\}} \frac{A_h(v_h, v_h)}{\|v_{\mathcal{T}}\|_{L^2}^2}
	\leq \max_{v_h\in \mathbb{E}_h\setminus\{0\}} \frac{A_h(v_h,v_h)}{\|v_{\mathcal{T}}\|_{L^2}^2}.
\end{align*}
The maximizer corresponds to $v_h = \mathcal{I} \phi$ for some $\phi\in \mathbb{E}$ with $\|\phi\|_{L^2}=1$.
The definition~\eqref{eq:modified_bilinear_form} of $A_h$ and the identity~\cref{eq:keyidentity} give
\begin{align}\label{eq:proof_majorization}
\begin{split}
\lambda_{h,m} \|\Pi_{\mathcal{T}}^{k+1}\phi\|_{L^2}^2
&\leq  \|\nabla_\pw \mathcal{G} \phi\|_{L^2}^2 + s_h(\mathcal{I} \phi, \mathcal{I} \phi)\\
&\qquad +  (V\Pi_{\mathcal{T}}^{k+1}\phi,\Pi_{\mathcal{T}}^{k+1}\phi)_{L^2}
- \alpha \|h_{\mathcal{T}}Q(V \Pi_{\mathcal{T}}^{k+1}\phi) \|_{L^2}^2.
\end{split}
\end{align}
For the first two terms on the right-hand side of~\eqref{eq:proof_majorization}, the Pythagoras identity $\|\nabla_\pw \mathcal{G} \phi\|_{L^2}^2 = \|\nabla \phi\|_{L^2}^2 - \|\nabla_\pw(\phi - \mathcal{G} \phi)\|_{L^2}^2$ and the stabilization bound~\eqref{eq:proof_stabilization} show
\begin{equation}\label{eq:proof_kinetic}
\|\nabla_\pw \mathcal{G} \phi\|_{L^2}^2 + s_h(\mathcal{I} \phi, \mathcal{I} \phi)
\leq \|\nabla \phi\|_{L^2}^2 + (\sigma\Lambda -1) \|\nabla_\pw(\phi - \mathcal{G}\phi)\|_{L^2}^2.
\end{equation}
To estimate the potential terms in~\eqref{eq:proof_majorization}, we use for $V_{\pm}=(\pm V)_+\geq0$ the identity
\begin{align*}
\|V_{\pm}^{1/2}\Pi_{\mathcal{T}}^{k+1}\phi\|_{L^2}^2 = \|V_{\pm}^{1/2}\phi\|_{L^2}^2 - 2(Q(V_{\pm}\Pi_{\mathcal{T}}^{k+1}\phi),Q\phi)_{L^2} - \|V_{\pm}^{1/2}Q\phi\|_{L^2}^2,
\end{align*}
which follows from the binomial formula and the $L^2$ orthogonality of $Q$.
Subtracting the identity for $V_-$ from that for $V_+$ gives the mixed term $-2(Q(V\Pi_{\mathcal T}^{k+1}\phi),Q\phi)_{L^2}$.
A Young inequality with weight $\alpha>0$ for this term and $\|V_+^{1/2}Q\phi\|_{L^2}^2\geq0$ imply
\begin{align*}%
	(V \Pi_{\mathcal{T}}^{k+1}\phi, \Pi_{\mathcal{T}}^{k+1}\phi)_{L^2} 
&- \alpha \|h_{\mathcal{T}} Q(V \Pi_{\mathcal{T}}^{k+1}\phi) \|_{L^2}^2\\
&\leq (V\phi,\phi)_{L^2} + \alpha^{-1}\|h_{\mathcal{T}}^{-1}Q\phi\|_{L^2}^2 +\|V_{-}^{1/2}Q\phi\|_{L^2}^2.
\end{align*}
Returning to~\eqref{eq:proof_majorization}, we insert~\eqref{eq:proof_kinetic} and the last estimate, and use the Pythagoras identity $\|\Pi_{\mathcal{T}}^{k+1}\phi\|_{L^2}^2 = 1 - \|Q\phi\|_{L^2}^2$ to obtain
\begin{align}\nonumber
	\lambda_{h,m} \leq \mathcal{E}(\phi)
		&+ (\sigma\Lambda -1) \|\nabla_\pw(\phi- \mathcal{G}\phi)\|_{L^2}^2\\
		&+ \alpha^{-1}\|h_{\mathcal{T}}^{-1}Q\phi\|_{L^2}^2+\lambda_{h,m}\|Q\phi\|_{L^2}^2+\|V_{-}^{1/2}Q\phi\|_{L^2}^2.\label{eq:proof_combined_rhs}
\end{align}
The Poincaré inequality~\eqref{eqn:CP} provides
\begin{align*}
	\alpha^{-1}\|h_{\mathcal{T}}^{-1}Q\phi\|_{L^2}^2
	+ \lambda_{h,m}\|Q\phi\|_{L^2}^2
	&\leq C_P^2(\alpha^{-1}+h_{\max}^2(\lambda_{h,m})_+)\|\nabla_\pw(\phi- \mathcal{G}\phi)\|_{L^2}^2.
\end{align*}
Using this and $\|V_-^{1/2}Q\phi\|_{L^2}\leq\delta_h\|\nabla_\pw(\phi-\mathcal G\phi)\|_{L^2}$ from~\eqref{eqn:V_negative_control} in~\eqref{eq:proof_combined_rhs} reveals
\begin{align*}%
	\lambda_{h,m} 
	  &\leq \mathcal{E}(\phi)
		+ \gamma(\lambda_{h,m}) \|\nabla_\pw(\phi- \mathcal{G}\phi)\|_{L^2}^2.
\end{align*}
This, $\mathcal{E}(\phi)\leq \lambda_m$,~\eqref{eqn:bestapprox_bound}, and $\lambda_m>0$ from~\eqref{eqn:Garding} conclude the proof.
\end{proof}

\begin{remark}[Choice of parameters]\label{rem:parameters}
	For practical applications, we balance the stabilization and discretization contributions by selecting $\sigma = \Lambda^{-1}/3$ and $\alpha = 3 C_P^2$, so that $\sigma\Lambda = C_P^2\alpha^{-1} = 1/3$. This particular choice yields
	\[
	\gamma(\lambda_{h,m}) = \bigl( h_{\max}^2 C_P^2\, (\lambda_{h,m})_{+} +\delta_h^2 - 1/3 \bigr)_+.
	\]
	Thus $\gamma(\lambda_{h,m})=0$ if and only if $h_{\max}^2C_P^2(\lambda_{h,m})_++\delta_h^2\leq1/3$.
For each fixed $m$ with bounded $(\lambda_{h,m})_+$, this condition holds on sufficiently fine meshes when $\delta_h=0$ for $V_-=0$ or when $\delta_h\in\mathcal O(h_{\max}^{1-n/(2p)})$ is chosen as in \Cref{lem:negative_potential_Lp}.
\end{remark}

\begin{remark}[Spectral shift]\label{rem:shift}
	One may always shift the spectrum by $M>0$ to make the $L^p$ norm of the negative potential part $V_-$ arbitrarily small, and thus obtain a coercive Schrödinger form by \Cref{lem:Garding_Lp}.
	Indeed, this shift amounts to replacing $V$ in~\eqref{eq:eigen_problem} by $\widehat V=V+M$ with negative part $\widehat V_-=(V_--M)_+$, satisfying
	\begin{align*}
		\|(V_--M)_+\|_{L^p}\to 0\qquad\text{as}\quad M\to \infty.
	\end{align*}
	Since $Q(Mv_{\mathcal T})=0$, the quadratic correction term in~\eqref{eq:modified_bilinear_form} is invariant under such a constant perturbation and the discrete spectrum is also exactly shifted by $M>0$.
	Hence \Cref{theo:main_lower_bound} immediately generalizes to
	\begin{equation*}
		\widehat{\mathrm{GLB}}(m) \coloneqq \frac{\widehat c_{A}(\lambda_{h,m}+M)_+}{\widehat c_{ A}+\widehat\gamma(\lambda_{h,m}+M)}-M \leq \lambda_m,
	\end{equation*}
	with coercivity constant $\widehat c_{A}>0$ and $\widehat\delta_h$ within $\widehat \gamma$ for the shifted potential $\widehat V$.

	Moreover, a spectral shift does not change the eigenfunctions.
	Hence the errors in the eigenvalues, eigenfunctions, and residuals are unaffected by this shift, and the outcome of the convergence analysis in \Cref{sec:apriori,sec:a-posteriori} transfers to general potentials $V$, with modified constants depending also on the coercivity shift $M$.
\end{remark}

\section{A priori error analysis}
\label{sec:apriori}

This section shows that the negative correction term in \cref{eq:modified_bilinear_form} is of higher order, so that the Babu\v{s}ka--Osborn theory~\cite{BabuskaOsborn1991} yields the same convergence rates as established for the HHO eigenvalue approximations in~\cite{CaloCicuttinErn2019}.
Throughout this and the following section, $A \lesssim B$ abbreviates $A \leq CB$ with a generic constant $C>0$ that is independent of the mesh size but may depend on $\Omega$, the shape regularity, $k$, $\sigma$, $\alpha$, $V$, the elliptic regularity constant, and the eigenvalue $\lambda$.

\subsection{Preparations}%
\label{sub:Higher-order potential correction}
To simplify the error analysis, oscillations of the potential contribution $f=Vv$ for conforming functions $v\in H^1_0(\Omega)$ are measured by
\begin{align}
	\mathrm{osc}_{k+1}(f,\mathcal{T})
	\coloneqq
	\|h_{\mathcal{T}}Qf\|_{L^2}
	\qquad\text{for all }f\in L^2(\Omega).
	\label{eq:oscillation}
\end{align}
We therefore assume that the potential $V$ is admissible in the following sense.
\begin{definition}[Admissible potential]\label{def:potential}
    A function $V \in L^p(\Omega)$ with $p > n/2$ and $p \ge 2$ is called an \emph{admissible potential} if multiplication by $V$ defines a bounded operator from $H_0^1(\Omega)$ to $L^2(\Omega)$, i.e., if there exists a constant $C_V > 0$ such that
    \begin{equation}\label{eq:continuous-multiplier}
        \|Vv\|_{L^2} \le C_V \|\nabla v\|_{L^2} \qquad\text{for all }v \in H_0^1(\Omega).
    \end{equation}
\end{definition}
H\"older's inequality and the Sobolev embedding verify~\eqref{eq:continuous-multiplier} if $p\geq n$ for $n\geq3$ or $p>2$ for $n=2$. 
Hardy's inequality further provides~\eqref{eq:continuous-multiplier} for the critical point singularity $V(x)=\pm|x-x_0|^{-1}\notin L^n(\Omega)$ at $x_0\in\Omega$ for $n\geq3$.
In particular, the attractive and repulsive Coulomb potentials are admissible in three dimensions.
We recall the following standard estimates without proof.
\begin{lemma}[Potential estimates]\label{lem:potential_estimates}
	Any $v\in H^1_0(\Omega)$ and any $v_T\in H^1(T)$ with $(v_T,1)_{L^2(T)}=0$ for $T\in\mathcal{T}$ satisfy
\begin{align}\label{eqn:multiplier_12}
    \||V|^{1/2}v\|_{L^2} 
	&\lesssim \|V\|_{L^p}^{1/2} \|\nabla v\|_{L^2},\\
	\label{eq:local_pot_est}
    \||V|^{1/2}v_T\|_{L^2(T)}
	&\lesssim h_T^{1-{n}/(2p)}
    \|V\|_{L^p(T)}^{1/2}
    \|\nabla v_T\|_{L^2(T)}.
\end{align}
\end{lemma}

The analysis relies on the smoothing operator of~\cite{EZ:QuasioptimalVariantHybrid2020}.
\begin{lemma}[Smoothing operator {\cite[Lem.~5.2]{CGT:AdaptiveHybridHighorder2024}}]\label{lem:enrichment}
There exist $k_{\mathcal J}\in\mathbb N$ and a linear right-inverse
$\mathcal J\colon V_h\to H_0^1(\Omega)\cap\mathbb P^{k_{\mathcal J}}(\mathcal{T})$
of the HHO interpolation $\mathcal I$, i.e.,
$\mathcal I\circ\mathcal J=\mathrm{id}$, such that the following holds.
	Any $v_h \in V_h$ and $v\in H^1_0(\Omega)$ satisfy
	\begin{align}
		\|\nabla \mathcal Jv_h\|_{L^2}
		&\lesssim \|v_h\|_{a_h},\label{eq:Jstability}\\
		\|\nabla(\mathcal J\mathcal Iv-v)\|_{L^2}
		&\lesssim \|\nabla_\mathrm{pw}(v-\mathcal{G} v)\|_{L^2}\label{eq:Japprox}
	\end{align}
	with hidden constants depending only on $\Omega$, the shape regularity, $k$, and~$\sigma$.
\end{lemma}

The next result shows that the discrete potential correction is of higher order.
\begin{lemma}[Correction term estimate]\label{lem:higher-order-defect}
    There exists $C_{\mathrm{cor}} > 0$, depending only on $\Omega$, the shape regularity, $k$, $\sigma$, and $V$, such that any $v_h \in V_h$ satisfies
    \begin{equation*}
        \|h_{\mathcal{T}} Q(V v_{\mathcal{T}})\|_{L^2}
        \le \|h_{\mathcal{T}} V v_{\mathcal{T}}\|_{L^2}
        \le C_{\mathrm{cor}} h_{\max}^{\min\{1,\, 2 - n/p\}} \|v_h\|_{a_h}.
    \end{equation*}
\end{lemma}
\begin{proof}
The first inequality follows from the $L^2$-stability of $Q$. For the second, the triangle inequality with the smoothing $\mathcal{J}v_h\in H_0^1(\Omega)$ of $v_h$ from \cref{lem:enrichment} yields
\begin{equation*}%
  \|h_{\mathcal{T}} V v_{\mathcal{T}}\|_{L^2}
  \le \|h_{\mathcal{T}} V \mathcal{J} v_h\|_{L^2}
  + \|h_{\mathcal{T}} V (\mathcal{J} v_h - v_{\mathcal{T}})\|_{L^2}.
\end{equation*}
To estimate the first term on the right-hand side, we apply the multiplier bound~\eqref{eq:continuous-multiplier} together with the stability estimate~\cref{eq:Jstability} for $\mathcal{J}$ to obtain
\begin{equation*}
  \|h_{\mathcal{T}} V \mathcal{J} v_h\|_{L^2}
  \le C_V h_{\max} \|\nabla (\mathcal{J} v_h)\|_{L^2}
  \lesssim h_{\max} \|v_h\|_{a_h}.
\end{equation*}
For the second term, $\mathcal I\circ\mathcal J=\mathrm{id}$ implies $\Pi_{\mathcal{T}}^{k+1}\mathcal{J}v_h=v_{\mathcal{T}}$ and hence $\mathcal{J}v_h-v_{\mathcal{T}}=Q\mathcal{J}v_h$. Hölder's inequality, an inverse estimate for the piecewise polynomial $Q\mathcal{J}v_h$, and~\eqref{eqn:CP} with the best approximation~\eqref{eq:bestapproximation} yield, for any $T\in\mathcal T$,
\begin{align*}
  h_T \|V Q\mathcal{J}v_h\|_{L^2(T)}
  &\le h_T \|V\|_{L^p(T)} \|Q\mathcal{J}v_h\|_{L^{2p/(p-2)}(T)}\\
  &\lesssim h_T^{1-n/p} \|V\|_{L^p(T)} \|Q\mathcal{J}v_h\|_{L^2(T)}\\
  &\lesssim h_T^{2-n/p} \|V\|_{L^p(T)} \|\nabla \mathcal{J}v_h\|_{L^2(T)}.
\end{align*}
Summing over all $T \in \mathcal{T}$ and applying the stability estimate~\cref{eq:Jstability} results in
\begin{equation*}
  \|h_{\mathcal{T}}V(\mathcal{J}v_h-v_{\mathcal{T}})\|_{L^2}
  \lesssim
  h_{\max}^{2-n/p}\|v_h\|_{a_h}.
\end{equation*}
Combining the above estimates completes the proof.
\end{proof}

Similar arguments bound the potential term $VQv$ by the elliptic projection error.
\begin{lemma}[Weighted projection error]\label{lem:local_pro_err}
    For any $v \in H_0^1(\Omega)$, we have
    \begin{equation*}
        \|h_{\mathcal{T}}VQ v\|_{L^2}
        \lesssim h_{\max}^{2-n/p} \|\nabla_{\mathrm{pw}}(v - \mathcal{G} v)\|_{L^2}.
    \end{equation*}
\end{lemma}

\begin{proof}
Since $Qv = v - \Pi_{\mathcal{T}}^{k+1} v$ and $\Pi_{\mathcal{T}}^{k+1}\mathcal{G}v=\mathcal{G}v$, the local $H^1$-stability of $\Pi_T^{k+1}$, as in \Cref{cor:H1_stability}, yields the auxiliary estimate
    \begin{align}\label{eq:Qh_H1_bound}
    \begin{split}
        \|\nabla_{\mathrm{pw}} Q v\|_{L^2}
        &\le \|\nabla_{\mathrm{pw}} (v - \mathcal{G} v)\|_{L^2} + \|\nabla_{\mathrm{pw}} \Pi_{\mathcal{T}}^{k+1} (\mathcal{G} v - v)\|_{L^2} \\
        &\lesssim \|\nabla_{\mathrm{pw}} (v - \mathcal{G} v)\|_{L^2}.
    \end{split}
    \end{align}
We distinguish two cases based on the integrability of the potential $V \in L^p(\Omega)$.
In the first case $p > n$, $Qv$ has vanishing integral mean on each $T\in\mathcal{T}$, so that H\"older's inequality and a scaled Sobolev--Poincar\'e inequality yield
\begin{equation*}
        h_T\|VQ v\|_{L^2(T)}
        \le h_T\|V\|_{L^p(T)} \|Q v\|_{L^{2p/(p-2)}(T)}
        \lesssim h_T^{2-n/p} \|V\|_{L^p(T)} \|\nabla Q v\|_{L^2(T)}.
\end{equation*}
Summing over all $T\in\mathcal{T}$ and using \eqref{eq:Qh_H1_bound} yields the claim in this case.

In the second case $n/2 < p \le n$, adding and subtracting the interpolant $\mathcal{J}\mathcal{I} v$ from \cref{lem:enrichment} and applying the triangle inequality yield
    \begin{equation*}
        \|h_{\mathcal{T}}VQ v\|_{L^2}
        \leq \|h_{\mathcal{T}}V(v - \mathcal{J}\mathcal{I} v)\|_{L^2} + \|h_{\mathcal{T}}V(\mathcal{J}\mathcal{I} v - \Pi_{\mathcal{T}}^{k+1} v)\|_{L^2}.
    \end{equation*}
    For the first term, \eqref{eq:continuous-multiplier}, \eqref{eq:Japprox}, and $h_{\max}\lesssim h_{\max}^{2-n/p}$ for $p\leq n$ show
    \begin{equation*}
        \|h_{\mathcal{T}}V(v - \mathcal{J}\mathcal{I} v)\|_{L^2} \le C_V h_{\max}\|\nabla(v - \mathcal{J}\mathcal{I} v)\|_{L^2} \lesssim h_{\max}^{2-n/p}\|\nabla_{\mathrm{pw}}(v - \mathcal{G} v)\|_{L^2}.
    \end{equation*}
    For the second term, since $\mathcal{J}\mathcal{I} v - \Pi_{\mathcal{T}}^{k+1}v = Q(\mathcal{J}\mathcal{I} v)$ is piecewise polynomial, H\"older's inequality, an inverse estimate, and~\eqref{eqn:CP} with~\eqref{eq:bestapproximation} give, for any $T\in\mathcal{T}$,
    \begin{align*}
        h_T\|V(\mathcal{J}\mathcal{I} v - \Pi_{\mathcal{T}}^{k+1} v)\|_{L^2(T)}
        &\lesssim h_T^{1-n/p}\|V\|_{L^p(T)}\|Q(\mathcal{J}\mathcal{I} v)\|_{L^2(T)}\\
        &\lesssim h_T^{2-n/p}\|V\|_{L^p(T)} \|\nabla(\mathcal{J}\mathcal{I} v - \Pi_{T}^{k+1} v)\|_{L^2(T)}.
    \end{align*}
    Summing over all $T\in\mathcal{T}$ and adding and subtracting $v$, we obtain
    \begin{equation*}
        \|h_{\mathcal{T}}V(\mathcal{J}\mathcal{I} v - \Pi_{\mathcal{T}}^{k+1} v)\|_{L^2}
        \lesssim h_{\max}^{2-n/p} \big(\|\nabla(\mathcal{J}\mathcal{I} v - v)\|_{L^2} + \|\nabla_{\mathrm{pw}}Q v\|_{L^2}\big).
    \end{equation*}
    The first term in the parentheses is bounded via~\cref{eq:Japprox} and the second term via~\eqref{eq:Qh_H1_bound}. Combining the previous bounds completes the proof.
\end{proof}
The coercivity of $A(\bullet,\bullet)$ is inherited by its discrete counterpart, also in a slightly strengthened form required for the source problem analysis in \Cref{lem:err_est:source} below.
\begin{lemma}[Discrete coercivity]\label{lem:shifted-coercivity}
    There exist mesh-independent constants $c_0,h_0>0$ such that $h_{\max}\leq h_0$ implies
    \begin{equation*}%
        A_h(v_h,v_h)
        -\frac12|v_h|_{s_h}^2
        \geq c_0\|v_h\|_{a_h}^2
        \qquad\text{for all }v_h\in V_h.
    \end{equation*}
\end{lemma}

\begin{proof}
	We argue by contradiction and discrete compactness:
	Assume that there exist a sequence of meshes indexed by $h$ with $h_{\max}\to0$ and functions $v_h\in V_h$ with
	\begin{align*}
	    \|v_h\|_{a_h}^2=1
		\qquad\text{and}\qquad
	    \limsup_{h_{\max}\to0}
	    \left(A_h(v_h,v_h)-\frac12|v_h|_{s_h}^2\right)\leq0.
	\end{align*}
	Stability of $\mathcal{J}$ from~\eqref{eq:Jstability} implies that $w_h\coloneqq\mathcal Jv_h$ is bounded in $H^1_0(\Omega)$.
	Weak compactness in $H^1_0(\Omega)$ and the Rellich embedding provide a subsequence (not relabelled) and a weak limit $w\in H^1_0(\Omega)$ with $w_h\rightharpoonup w$ in $H_0^1(\Omega)$ and $w_h\to w$ in $L^2(\Omega)$. 
	Throughout this proof, convergence is understood as $h_{\max}\to0$.

	The identities $\mathcal I\mathcal J=\mathrm{id}$ and $\mathcal R\mathcal I=\mathcal G$ reveal $\mathcal Rv_h=\mathcal Gw_h$.
	Hence $\nabla_{\mathrm{pw}}\mathcal Rv_h-\nabla w_h$ is orthogonal to all elementwise constant vector fields.
	Since this difference is bounded in $L^2$ and $\Pi_{\mathcal T}^0\psi\to\psi$ for every fixed $\psi\in L^2(\Omega;\mathbb R^n)$, it follows that $\nabla_{\mathrm{pw}}\mathcal Rv_h\rightharpoonup\nabla w$ in $L^2$.
	Weak lower semicontinuity of the norm reveals $\|\nabla w\|_{L^2}^2\leq1$.

	Recall $Qw_h= w_h - v_{\mathcal{T}}$ from $\mathcal I\circ\mathcal J=\mathrm{id}$, which satisfies $\|h_{\mathcal T}^{-1}Qw_h\|_{L^2}\lesssim1$ by~\eqref{eq:bestapproximation} and~\eqref{eqn:CP}.
	Hence the multiplier estimate~\eqref{eq:continuous-multiplier} and \Cref{lem:higher-order-defect} provide
	\begin{align*}
	    &\bigl|(Vv_{\mathcal T},v_{\mathcal T})_{L^2}
	           -(Vw_h,w_h)_{L^2}\bigr|
			   \leq
	    \left(\|h_{\mathcal T}Vv_{\mathcal T}\|_{L^2}
	          +h_{\max}\|Vw_h\|_{L^2}\right)
	    \|h_{\mathcal T}^{-1}Qw_h\|_{L^2}
	    \to0.
	\end{align*}
	Moreover, the multiplier estimate and strong $L^2$ convergence lead to
\begin{align*}
    \bigl|(Vw_h,w_h)_{L^2}-(Vw,w)_{L^2}\bigr|
    \leq C_V\|\nabla(w_h+w)\|_{L^2}\|w_h-w\|_{L^2}
    \to0.
\end{align*}
Thus $(Vv_{\mathcal T},v_{\mathcal T})_{L^2}\to(Vw,w)_{L^2}$.
The potential correction term also tends to zero by \Cref{lem:higher-order-defect}.
This and the weak lower semicontinuity of norms with $1=\|v_h\|_{a_h}^2=\|\nabla_\pw\mathcal{R}v_h\|_{L^2}^2+|v_h|_{s_h}^2$ by definition imply
\begin{align*}
    0
    &\geq \liminf_{h_{\max}\to0}
    \left(A_h(v_h,v_h)-\frac12|v_h|_{s_h}^2\right)
    \geq \frac12+\frac12\|\nabla w\|_{L^2}^2+(Vw,w)_{L^2}\\
    &=A(w,w)+\frac12\left(1-\|\nabla w\|_{L^2}^2\right)
    \geq\min\{c_A,1/2\}>0,
\end{align*}
with continuous coercivity~\eqref{eqn:Garding} and $0\leq\|\nabla w\|_{L^2}^2\leq1$ in the last step.
This contradiction concludes the proof.
\end{proof}

\subsection{Spectral correctness}\label{sub:spectral}
By the Babu\v{s}ka--Osborn theory, spectral correctness follows from the convergence of the solution operators of the continuous and discrete source problems in the operator norm. 
Throughout the convergence analysis in \Cref{sec:apriori,sec:a-posteriori}, let $h_{\max}\leq h_0$ with $h_0$ from \Cref{lem:shifted-coercivity}.
The coercivity of $A(\bullet,\bullet)$ and $A_h(\bullet,\bullet)$ ensures that the continuous and discrete source problems are well posed.

The continuous solution operator
$T\colon L^2(\Omega)\to L^2(\Omega)$ maps any $f\in L^2(\Omega)$ to the unique solution $Tf\in H_0^1(\Omega)$ of
\begin{equation}\label{eq:source_pro}
    A(Tf,v)
    =b(f,v)
    \qquad\text{for all }v\in H_0^1(\Omega).
\end{equation}
If $(\lambda,u)$ is an eigenpair of \eqref{eq:weak_form}, then $(\lambda^{-1},u)$ is an eigenpair of $T$, and we denote the eigenspace of $\lambda$ by $E_\lambda\coloneqq\ker(\lambda^{-1}I-T)$.
For $h_{\max}\leq h_0$, the discrete solution operator $S_h\colon L^2(\Omega)\to V_h$ maps any $f\in L^2(\Omega)$ to the unique solution $S_hf\in V_h$ of
\begin{equation}\label{eq:dis_source_pro}
    A_h(S_hf,v_h)
	=b(f,v_{\mathcal{T}})
    \qquad\text{for all }v_h\in V_h.
\end{equation}
The cell component of $S_hf$ defines ${T_h\colon L^2(\Omega)\to\mathbb{P}^{k+1}(\mathcal{T})}$ by $T_hf\coloneqq(S_hf)_{\mathcal{T}}$.
If $(\lambda_h,u_h)$ is an eigenpair of \eqref{eq:discrete_eigenproblem_modif}, then $(\lambda_h^{-1},u_{\mathcal{T}})$ is an eigenpair of $T_h$.

To obtain quantitative convergence estimates, assume that there exists a regularity index $s\in(0,1]$ such that all $0\leq t\leq s$ satisfy
\begin{equation} \label{eq:ellipticregularity}
\|Tf\|_{H^{1+t}} \lesssim \|f\|_{H^{t-1}}
\qquad\text{for all }f\in L^2(\Omega).
\end{equation}
Such a result is classical for the Laplacian, see, e.g., \cite[Thm.~5.5]{AmroucheMoussaoui2025}, and transfers to the Schr\"odinger operator using the multiplier bound~\eqref{eq:continuous-multiplier}.
The following lemma relies on the source problem estimates from \cref{sec:apriori:A}.

\begin{lemma}[Spectral correctness]\label{lem:norm_conver}
    Let $h_{\max}$ be sufficiently small. Then
    \begin{equation}
    \label{eq:approxopnorm}
        \|T - T_h\|_{\mathcal{L}(L^2; L^2)} \lesssim h_{\max}^s.
    \end{equation}
    For any eigenvalue $\lambda=\lambda_m=\cdots=\lambda_{m+r-1}$ of \eqref{eq:weak_form} with multiplicity $r$ and $\mu\coloneqq\lambda^{-1}$, the eigenvalues $\mu_{h,j}\coloneqq\lambda_{h,j}^{-1}$ of $T_h$ satisfy
    \begin{equation}
    \label{eq:approxev}
        \max_{m \le j \le m+r-1} |\mu - \mu_{h,j}| \lesssim \sup_{\phi,\psi \in E_{\lambda} \setminus \{0\}} \frac{|((T - T_h)\phi, \psi)_{L^2}|}{\|\phi\|_{L^2} \|\psi\|_{L^2}} + \|(T - T_h)|_{E_{\lambda}}\|_{\mathcal{L}(E_{\lambda}; L^2)}^2.
    \end{equation}
    Moreover, if $u_{\mathcal{T},j} \in \mathbb{P}^{k+1}(\mathcal{T})$ with $\|u_{\mathcal{T},j}\|_{L^2} = 1$ is an eigenvector associated with $\mu_{h,j}$ for $m\le j\le m+r-1$, then there exists a unit eigenvector $u \in E_{\lambda}$ such that
    \begin{equation}\label{eq:approxef}
        \|u - u_{\mathcal{T},j}\|_{L^2} \lesssim \|(T - T_h)|_{E_{\lambda}}\|_{\mathcal{L}(E_{\lambda}; L^2)}.
    \end{equation}
\end{lemma}
\begin{proof}
    Let $f \in L^2(\Omega)$. The triangle inequality, \eqref{eqn:CP}, \cref{lem:L2err_est:source}, and~\eqref{eq:bestapproximation} yield
    \begin{align*}
        \|Tf - T_h f\|_{L^2}
        &\le \|Q Tf\|_{L^2} + \|\Pi_{\mathcal{T}}^{k+1} Tf - T_h f\|_{L^2} \\
        &\lesssim h_{\max} \|\nabla Tf\|_{L^2} + h_{\max}^s \big( \|\nabla Tf\|_{L^2} + \mathrm{osc}_{k+1}(f, \mathcal{T})\\
        &\qquad + \mathrm{osc}_{k+1}(V(Tf), \mathcal{T}) \big) \\
        &\lesssim h_{\max}^s \big( \|\nabla Tf\|_{L^2} + h_{\max}\|f\|_{L^2} \big),
    \end{align*}
    where the last step uses the $L^2$-stability of $Q$ and the multiplier estimate~\eqref{eq:continuous-multiplier}. Since $\|\nabla Tf\|_{L^2}\lesssim\|f\|_{L^2}$, this establishes \cref{eq:approxopnorm}.
    The operator $T$ is compact by the Rellich embedding, $T_h$ has finite rank, and both operators are self-adjoint and positive semidefinite by the symmetry and coercivity of $A$ and $A_h$. Hence, the max--min principle and \cref{eq:approxopnorm} yield $|\lambda_m^{-1}-\lambda_{h,m}^{-1}|\le\|T-T_h\|_{\mathcal{L}(L^2;L^2)}\lesssim h_{\max}^s$ for all $1\leq m\leq N$, which justifies the indexing, and the remaining estimates \cref{eq:approxef,eq:approxev} follow from the Babu\v{s}ka--Osborn theory \cite[Thm.~7.3 \& 7.4]{BabuskaOsborn1991}. Since $\|u_{\mathcal{T},j}\|_{L^2}=1$, normalizing $u$ in \cref{eq:approxef} only changes the hidden constant.
\end{proof}

\subsection{A priori error estimates for eigenpairs}\label{sub:apriori_eigenpairs}
The main result shows that the eigenvalue errors converge at twice the rate of the eigenvector errors in the energy norm.
Recall the elliptic regularity index $s\in(0,1]$ from~\cref{eq:ellipticregularity} and set
\begin{equation}
\label{eq:sv}
    s_V \coloneqq \min\left\{s,1-\frac{n}{2p}\right\}>0,
\end{equation}
which quantifies the gain of the $|V|^{1/2}$-weighted $L^2$ error over the energy error. The best-approximation error and potential oscillation of $v\in H^1_0(\Omega)$ are combined in
\begin{equation}\label{eq:Err}
    \operatorname{Err}(v)\coloneqq \|\nabla_{\mathrm{pw}}(v - \mathcal{G} v)\|_{L^2}^2 + \mathrm{osc}_{k+1}^2(Vv, \mathcal{T}).
\end{equation}

\begin{theorem}[A priori error estimates for eigenpairs]\label{thm:a_priori}
    Let $\lambda=\lambda_m=\cdots=\lambda_{m+r-1}$ be an eigenvalue of \eqref{eq:weak_form} with multiplicity $r$.
For sufficiently small $h_{\max}$, the discrete eigenvalues $\lambda_{h,m}, \dots, \lambda_{h,m+r-1}$ of \eqref{eq:discrete_eigenproblem_modif} satisfy
    \begin{equation*}
        \max_{m \le j \le m+r-1} |\lambda - \lambda_{h,j}| \lesssim \sup_{\substack{v \in E_\lambda \\ \|v\|_{L^2} = 1}} \operatorname{Err}(v).
    \end{equation*}
    If $u_{h,j} = (u_{\mathcal{T},j}, u_{\mathcal{F},j}) \in V_h$ with $\|u_{\mathcal{T},j}\|_{L^2} = 1$ is an eigenvector associated with $\lambda_{h,j}$ for some $m \le j \le m+r-1$, then there exists a unit eigenvector $u \in E_{\lambda}$ such that
    \begin{equation*}
        \|\mathcal{I} u - u_{h,j}\|_{a_h}^2
        + h_{\max}^{-2s_V} \||V|^{1/2}(u_{\mathcal{T},j} - u)\|_{L^2}^2
        + h_{\max}^{-2s} \|u_{\mathcal{T},j} - u\|_{L^2}^2 \lesssim \sup_{\substack{v \in E_\lambda \\ \|v\|_{L^2} = 1}} \operatorname{Err}(v).
    \end{equation*}
\end{theorem}

For sufficiently smooth eigenfunctions and potentials, standard approximation properties bound the right-hand sides in \cref{thm:a_priori} by $h_{\max}^{2(k+1)}$, so that the eigenvalue errors are of order $2(k+1)$ and the eigenvector errors in the energy norm of order $k+1$. This recovers the rates known from the HHO method~\cite{CaloCicuttinErn2019}.

\begin{proof}[Proof of \cref{thm:a_priori}]
We first bound the operator norm of the difference $T-T_h$ on the eigenspace $E_\lambda$. Any $v\in E_\lambda$ satisfies $\operatorname{osc}_{k+1}(\lambda v,\mathcal T)\lesssim h_{\max}^2\|\nabla_\mathrm{pw}(v-\mathcal Gv)\|_{L^2}$ and $\|QTv\|_{L^2}\lesssim h_{\max}\|\nabla_\mathrm{pw}(v-\mathcal Gv)\|_{L^2}$ by \cref{eqn:CP}. Hence the triangle inequality, which splits off $QTv$, and \cref{lem:L2err_est:source} for the remaining projected error yield the estimate
    \begin{align}\label{eq:opnormE}
    \|(T-T_h)|_{E_{\lambda}}\|_{\mathcal{L}(E_{\lambda};L^2)}
    \lesssim h_{\max}^s
    \sup_{\substack{v\in E_\lambda\\ \|v\|_{L^2}=1}}\operatorname{Err}(v)^{1/2}.
    \end{align}

    For fixed $j\in\{m,\dots,m+r-1\}$, we write $(\lambda_h,u_h)\coloneqq(\lambda_{h,j},u_{h,j})$ with $u_h = (u_{\mathcal{T}},u_{\mathcal{F}})$. Then, \cref{eq:approxef,eq:opnormE} give a unit eigenvector $u\in E_\lambda$ with
    \begin{align}
    \label{eq:eigenvalueL2conv}
        \|u - u_{\mathcal{T}}\|_{L^2}\lesssim h_{\max}^s\sup_{\substack{v\in E_\lambda\\ \|v\|_{L^2}=1}}\operatorname{Err}(v)^{1/2}.
    \end{align}

    Next, let $\phi,\psi\in E_\lambda$ with $\|\phi\|_{L^2}=\|\psi\|_{L^2}=1$. By the definitions of $T_h$ and $S_h$ in~\cref{eq:dis_source_pro} and the symmetry of $A_h$, the difference $T-T_h$ satisfies
    \begin{align}
        ((T-T_h)\phi,\psi)_{L^2}
        &= (T\phi,\psi)_{L^2} - A_h(S_h\phi,S_h\psi)\notag\\
        &= (T\phi,\psi)_{L^2}
        - A_h(\mathcal{I}T\phi,S_h\psi)
        + A_h(\mathcal{I}T\phi-S_h\phi,S_h\psi)\notag\\
        &= (T\phi,\psi)_{L^2}
        - (\psi,\Pi_{\mathcal{T}}^{k+1}T\phi)_{L^2}
        + A_h(\mathcal{I}T\phi-S_h\phi,S_h\psi).
        \label{eq:operator_diff}
    \end{align}
    Noting that $T\phi=\tfrac{1}{\lambda}\phi$, the $L^2$-orthogonality properties of $Q$, together with \cref{eqn:CP}, yield the following bound for the first two terms on the right-hand side of \cref{eq:operator_diff}:
    \begin{align*}
        (T\phi,\psi)_{L^2} - (\psi,\Pi_{\mathcal{T}}^{k+1}T\phi)_{L^2}
        &= \tfrac{1}{\lambda}(Q\phi,Q\psi)_{L^2}\\
        &\lesssim h_{\max}^2
        \|\nabla_{\mathrm{pw}}(\phi-\mathcal{G}\phi)\|_{L^2}
        \|\nabla_{\mathrm{pw}}(\psi-\mathcal{G}\psi)\|_{L^2}.
    \end{align*}
    The last term of \cref{eq:operator_diff} can be further split into three components:
    \begin{align}
    \begin{split}\label{eq:bh_expansion}
        A_h(\mathcal{I}T\phi - S_h\phi,S_h\psi) &= A_h(\mathcal{I}T\phi - S_h\phi,S_h\psi - \mathcal{I}T\psi)\\
        &\quad + \big( A(T\phi,T\psi) - A_h(S_h\phi,\mathcal{I}T\psi) \big) \\
        &\quad + \big( A_h(\mathcal{I}T\phi ,\mathcal{I}T\psi) - A(T\phi,T\psi) \big).
    \end{split}
    \end{align}
    Using the definition of $A_h$, together with \cref{lem:higher-order-defect,lem:err_est:source}, we obtain the following bound for the first term on the right-hand side of \cref{eq:bh_expansion}:
    \begin{equation*}
        A_h(\mathcal{I}T\phi - S_h\phi,S_h\psi - \mathcal{I}T\psi)
        \lesssim \operatorname{Err}(\phi)^{1/2}\operatorname{Err}(\psi)^{1/2},
    \end{equation*}
    where we used the relations $T\phi=\tfrac{1}{\lambda}\phi$ and $T\psi=\tfrac{1}{\lambda}\psi$. For the second term on the right-hand side of \cref{eq:bh_expansion}, the definitions of $T$ and $S_h$ in \cref{eq:source_pro,eq:dis_source_pro} show $A(T\phi,T\psi) - A_h(S_h\phi,\mathcal{I}T\psi) = (\phi,QT\psi)_{L^2} = \tfrac{1}{\lambda}(Q\phi,Q\psi)_{L^2}$, which is bounded as the first two terms in \cref{eq:operator_diff}.
    For the last term in \cref{eq:bh_expansion}, we write
    \begin{align}
        &A_h(\mathcal{I}T\phi,\mathcal{I}T\psi)
        - A(T\phi,T\psi)\notag\\
        \begin{split}
            &\quad = \tfrac{1}{\lambda^2} \Big(
            \big(
                (\nabla_\mathrm{pw}\mathcal{G}\phi,\nabla_\mathrm{pw}\mathcal{G}\psi)_{L^2}
                -(\nabla\phi,\nabla\psi)_{L^2}
            \big)\\
            &\qquad\qquad
            + \big(
                (V\Pi_{\mathcal{T}}^{k+1}\phi,\Pi_{\mathcal{T}}^{k+1}\psi)_{L^2}
                -(V\phi,\psi)_{L^2}
            \big)\\
            &\qquad\qquad
            + 
                s_h(\mathcal{I}\phi,\mathcal{I}\psi)
                -\alpha
                (h_{\mathcal{T}}Q(V\Pi_{\mathcal{T}}^{k+1}\phi),
                h_{\mathcal{T}}Q(V\Pi_{\mathcal{T}}^{k+1}\psi))_{L^2}
            \Big).
        \end{split}
        \label{eq:bh_b_diff}
    \end{align}
    We bound the four terms on the right-hand side of \cref{eq:bh_b_diff} separately. First, the orthogonality properties of the elliptic projection in \cref{eq:defGh} yield
    \begin{align*}
        (\nabla_\mathrm{pw}\mathcal{G}\phi, \nabla_\mathrm{pw}\mathcal{G}\psi)_{L^2}- (\nabla\phi,\nabla\psi)_{L^2}
        &= (\nabla_\mathrm{pw}(\mathcal{G}\phi-\phi),\nabla_\mathrm{pw}(\psi - \mathcal{G}\psi))_{L^2}\\
        &\leq \|\nabla_\mathrm{pw}(\mathcal{G}\phi-\phi)\|_{L^2}\|\nabla_\mathrm{pw}(\mathcal{G}\psi-\psi)\|_{L^2}.
    \end{align*}
    For the second term in \cref{eq:bh_b_diff}, the $L^2$-orthogonality properties of $Q$, the splitting $(VQ\phi,Q\psi)_{L^2}=(h_{\mathcal{T}}VQ\phi,h_{\mathcal{T}}^{-1}Q\psi)_{L^2}$ with \cref{lem:local_pro_err}, and \cref{eqn:CP} imply
    \begin{align*}
        &(V\Pi_{\mathcal{T}}^{k+1} \phi,\Pi_{\mathcal{T}}^{k+1}\psi)_{L^2} - (V\phi,\psi)_{L^2}\\
        &\quad = - (Q\phi,Q(V\psi))_{L^2} + (VQ\phi,Q\psi)_{L^2} - (Q(V\phi),Q\psi)_{L^2}\\
        &\qquad \lesssim \operatorname{Err}(\phi)^{1/2}\operatorname{Err}(\psi)^{1/2}.
    \end{align*}
    The stabilization term can be bounded using \cref{eq:proof_stabilization} as
    \begin{equation*}
        s_h(\mathcal{I}\phi,\mathcal{I}\psi)\leq|\mathcal{I}\phi|_{s_h}|\mathcal{I}\psi|_{s_h}\lesssim \|\nabla_\mathrm{pw}(\mathcal{G}\phi-\phi)\|_{L^2}\|\nabla_\mathrm{pw}(\mathcal{G}\psi-\psi)\|_{L^2}.
    \end{equation*}
    For the last term, $\Pi_{\mathcal{T}}^{k+1}\phi=\phi-Q\phi$, the $L^2$-stability of $Q$, and \cref{lem:local_pro_err} yield
    \begin{align*}
        &(h_{\mathcal{T}} Q (V \Pi_{\mathcal{T}}^{k+1}\phi), h_{\mathcal{T}} Q(V \Pi_{\mathcal{T}}^{k+1}\psi))_{L^2}\\
        &\quad \lesssim  \big(h_{\max}^{2-n/p}\|\nabla_\mathrm{pw}(\mathcal{G}\phi-\phi)\|_{L^2}
        +\mathrm{osc}_{k+1}(V\phi,\mathcal T) \big)\\
        &\qquad \times \big(h_{\max}^{2-n/p}\|\nabla_\mathrm{pw}(\mathcal{G}\psi-\psi)\|_{L^2}
        +\mathrm{osc}_{k+1}(V\psi,\mathcal T)\big).
    \end{align*}

    Combining the estimates for the terms in \eqref{eq:operator_diff}--\eqref{eq:bh_b_diff} yields
    \begin{equation*}
        ((T-T_h)\phi,\psi)_{L^2} \lesssim \operatorname{Err}(\phi)^{1/2}\operatorname{Err}(\psi)^{1/2}.
    \end{equation*}
Then, by \cref{lem:norm_conver} and \cref{eq:opnormE}, and noting that $\lambda=\mu^{-1}$ and $\lambda_{h,j}=(\mu_{h,j})^{-1}$ for $m\leq j\leq m+r-1$, we obtain the following eigenvalue estimate:    \begin{equation}\label{eq:eigenvalueconv}
        \max_{m \leq j \leq m+r-1}|\lambda - \lambda_{h,j}|\lesssim \sup_{\substack{v\in E_\lambda\\ \|v\|_{L^2}=1}}\operatorname{Err}(v).
    \end{equation}

    Next, we consider the estimate of $\|\mathcal{I}u-u_h\|_{a_h}$, where $u\in E_\lambda$ with $\|u\|_{L^2}=1$ is as given above. Abbreviating $e_h = \mathcal{I} u -u_h = (e_{\mathcal{T}},e_{\mathcal{F}})$ and noting that
    \begin{equation*}
        2\lambda_h = \lambda_h \|u-u_{\mathcal{T}}\|_{L^2}^2 +2\lambda_h(u,u_{\mathcal{T}})_{L^2} = \lambda_h \|u-u_{\mathcal{T}}\|_{L^2}^2 +2A_h(u_h,\mathcal{I}u),
    \end{equation*}
    we find, with $\lambda=A(u,u)$ and $\lambda_h=A_h(u_h,u_h)$, that
    \begin{align}
        A_h(e_h,e_h) &= \lambda_h \|u-u_{\mathcal{T}}\|_{L^2}^2 +\big(\lambda-\lambda_h\big)+A_h(\mathcal{I}u,\mathcal{I}u) - \lambda.\label{eq:estbhh}
    \end{align}
Using the coercivity of $A_h$ and combining \cref{eq:eigenvalueL2conv,eq:eigenvalueconv} to bound the first and second terms on the right-hand side of \cref{eq:estbhh}, while applying a similar argument to the expansion in \cref{eq:bh_b_diff} to bound the last term, we deduce
\begin{equation*}
        \|e_h\|_{a_h}^2\lesssim\sup_{\substack{v\in E_\lambda\\ \|v\|_{L^2}=1}}\operatorname{Err}(v).
    \end{equation*}

Finally, we prove the improved error estimate for $\||V|^{1/2}(u_{\mathcal{T}}-u)\|_{L^2}^2$. The triangle inequality and the splitting $u_{\mathcal{T}}-u=-Qu-(e_{\mathcal{T}}-\mathcal{J}e_h)-\mathcal{J}e_h$ yield
\begin{equation}
\label{eq:VpotimpL2}
        \||V|^{1/2}(u_{\mathcal{T}}-u)\|_{L^2}\leq \||V|^{1/2}Qu\|_{L^2} + \||V|^{1/2}(e_{\mathcal{T}}-\mathcal{J}e_h)\|_{L^2} + \||V|^{1/2}\mathcal{J}e_h\|_{L^2}.
    \end{equation}
For the first term on the right-hand side of~\cref{eq:VpotimpL2},  \cref{eq:local_pot_est,eq:Qh_H1_bound} yield
\begin{equation*}
    \||V|^{1/2}Qu\|_{L^2}
    \lesssim h_{\max}^{1-n/(2p)}
    \|\nabla_\mathrm{pw}(u-\mathcal{G}u)\|_{L^2},
\end{equation*}
while for the second term, the identity $e_{\mathcal{T}}-\mathcal{J}e_h=-Q\mathcal{J}e_h$ from $\mathcal I\circ\mathcal J=\mathrm{id}$, together with \cref{eq:local_pot_est}, \cref{eq:Qh_H1_bound} with \cref{eq:bestapproximation}, and \cref{eq:Jstability}, gives
\begin{equation*}
    \||V|^{1/2}(e_{\mathcal{T}}-\mathcal{J}e_h)\|_{L^2}
    \lesssim h_{\max}^{1-n/(2p)}
    \|\nabla\mathcal{J}e_h\|_{L^2}
    \lesssim h_{\max}^{1-n/(2p)}\|e_h\|_{a_h}.
\end{equation*}
For the last term in \cref{eq:VpotimpL2}, we employ an Aubin--Nitsche-type duality argument. Consider the dual problem $z \coloneqq T(|V|\mathcal{J}e_h) \in H_0^1(\Omega)$, with $T$ defined in \cref{eq:source_pro}. By the elliptic regularity estimate \cref{eq:ellipticregularity} with $t=s_V$, we have
\[
\|z\|_{H^{s_V+1}}
\lesssim \||V|\mathcal{J}e_h\|_{H^{s_V-1}}
\lesssim \||V|^{1/2}\mathcal{J}e_h\|_{L^2}.
\]
In the last step, we factor $|V|=|V|^{1/2}|V|^{1/2}$ and absorb one factor into the test function, applying the Sobolev inequality similarly as in \cref{eqn:multiplier_12}. The Sobolev embedding theorem ensures the continuous embedding $H^{1-s_V}(\Omega)\hookrightarrow L^{2p/(p-1)}(\Omega)$, since the condition $1-s_V\geq n/(2p)$ holds by \cref{eq:sv}.

Testing the dual problem with $\mathcal{J}e_h$ and inserting $\mathcal{R}e_h$ gives
\begin{align*}
    \||V|^{1/2}\mathcal{J}e_h\|_{L^2}^2
    &= (|V|\mathcal{J}e_h,\mathcal{J}e_h)_{L^2}
    = (\nabla z,\nabla\mathcal{J}e_h)_{L^2}
    +(Vz,\mathcal{J}e_h)_{L^2}\\
    &= (\nabla z,\nabla_\mathrm{pw}(\mathcal{J}e_h-\mathcal{R}e_h))_{L^2}
    + (\nabla z,\nabla_\mathrm{pw}\mathcal{R}e_h)_{L^2}
    +(Vz,\mathcal{J}e_h)_{L^2},
\end{align*}
which has the same structure as \cref{eq:rhsopL21} in the proof of \cref{lem:L2err_est:source}. We follow the arguments in that proof. Therein, \cref{lem:err_est:source} is replaced by the above bound for $\|e_h\|_{a_h}$ and $|u_h|_{s_h}\le|e_h|_{s_h}+|\mathcal{I}u|_{s_h}$, \cref{eq:estG} by $\|\nabla_\pw(z-\mathcal Gz)\|_{L^2}\lesssim h_{\max}^{s_V}\||V|^{1/2}\mathcal Je_h\|_{L^2}$, and the estimate in \cref{eq:rhsdefect} by
\begin{align*}
    &(\lambda u,\mathcal{J}\mathcal{I}z)_{L^2}
    -(\lambda_hu_{\mathcal{T}},\Pi_{\mathcal{T}}^{k+1}z)_{L^2}=(\lambda u-\lambda_hu_{\mathcal{T}},\mathcal{J}\mathcal{I}z)_{L^2}\\
    &\qquad =(\lambda-\lambda_h)(u,\mathcal{J}\mathcal{I}z)_{L^2}
    +\lambda_h(u-u_{\mathcal{T}},\mathcal{J}\mathcal{I}z)_{L^2}\\
    &\qquad \lesssim
    \bigl(|\lambda-\lambda_h|+\|u-u_{\mathcal{T}}\|_{L^2}\bigr)
    \|z\|_{H^1}.
\end{align*}
For all $v\in E_\lambda$ with $\|v\|_{L^2}=1$, \cref{eq:ellipticregularity} for $v=T(\lambda v)$ and \cref{eq:continuous-multiplier} show
\begin{align}\label{eqn:Err_bound}
	\operatorname{Err}(v)\lesssim h_{\max}^{2s}.
\end{align}
Together with \cref{eq:eigenvalueL2conv,eq:eigenvalueconv}, this yields for the last term in \cref{eq:VpotimpL2}
\begin{equation*}
    \||V|^{1/2}\mathcal{J}e_h\|_{L^2}
    \lesssim h_{\max}^{s_V}
    \sup_{\substack{v\in E_\lambda\\ \|v\|_{L^2}=1}}\operatorname{Err}(v)^{1/2}.
\end{equation*}
Combining \cref{eq:VpotimpL2} with these bounds and the bound for $\|e_h\|_{a_h}$ concludes the proof.
\end{proof}

\section{A posteriori error analysis}\label{sec:a-posteriori}

This section derives a reliable and efficient a posteriori error estimator for a simple eigenvalue in $n\in\{2,3\}$ dimensions. To this end, it adapts the stabilization-free a posteriori framework of~\cite{BCGT:StabilizationfreeHHOPosteriori2023,CGT:AdaptiveHybridHighorder2024} to spatially varying potentials and the additional correction term of the modified HHO scheme.
Let $(\lambda,u)\coloneqq(\lambda_m,u_m)$ denote a simple eigenpair of~\eqref{eq:weak_form} and $(\lambda_h,u_h)\coloneqq(\lambda_{h,m},u_{h,m})$ the corresponding discrete eigenpair of~\eqref{eq:discrete_eigenproblem_modif} with $(u,u_{\mathcal T})_{L^2}\geq 0$.
The estimator uses the discrete gradient $p_h\coloneqq \nabla_\pw\mathcal{R}u_h$ and reads $\eta^2(\mathcal T)\coloneqq \sum_{T\in\mathcal T} \eta^2(T)$ with local contributions
\begin{align}
	\eta^2(T)
	&\coloneqq
	|T|^{2/n}
	\|\operatorname{div} p_h+\lambda_hu_{\mathcal T}-Vu_{\mathcal T}\|_{L^2(T)}^2
	+|T|^{1/n}\sum_{F\in\mathcal{F}(T)}
	\|[p_h]_F\|_{L^2(F)}^2.
	\label{eq:local-estimator}
\end{align}
Here $[\bullet]_F$ abbreviates the full jump on interior faces $F\in\mathcal F(\Omega)$ and the tangential trace on boundary faces $F\in\mathcal F(\partial \Omega)$.
The framework of~\cite{BCGT:StabilizationfreeHHOPosteriori2023} relies on two equilibrium properties of $p_h$ for the Courant and Raviart--Thomas spaces
\begin{align*}
	S_0^1(\mathcal T)
	&\coloneqq \mathbb P^1(\mathcal T)\cap H_0^1(\Omega),\\
	\mathrm{RT}_0(\mathcal T)
	&\coloneqq
	\{q\in H(\operatorname{div},\Omega)
	\mid q|_T\in \mathbb P^0(T;\mathbb R^n)+x\mathbb P^0(T)
		\text{ for all }T\in\mathcal T\}.
\end{align*}

\begin{lemma}[Discrete equilibrium identities]\label{lem:a-posteriori-equilibrium}
	For any $v_C\in S_0^1(\mathcal T)$ and any $q_{\mathrm{RT}}\in\mathrm{RT}_0(\mathcal T)$ with $\operatorname{div}q_{\mathrm{RT}}=0$, the discrete gradient $p_h=\nabla_\pw \mathcal{R}u_h$ satisfies
	\begin{align}
		(p_h,\nabla v_C)_{L^2}
		&=((\lambda_h-V)u_{\mathcal T},v_C)_{L^2}
		+\alpha
		(h_{\mathcal{T}} Q(Vu_{\mathcal T}),h_{\mathcal{T}} Q(Vv_C))_{L^2},\label{eq:equilibrium-S}\\
		(p_h,q_{\mathrm{RT}})_{L^2}&=0.\label{eq:equilibrium-RT}
	\end{align}
\end{lemma}

\begin{proof}
	Polynomial consistency implies $\mathcal R\mathcal Iv_C=v_C$ and $s_h(u_h,\mathcal Iv_C)=0$.
	Hence $a_h(u_h,\mathcal Iv_C)=(p_h,\nabla v_C)_{L^2}$, and the discrete eigenvalue equation~\eqref{eq:discrete_eigenproblem_modif} proves~\eqref{eq:equilibrium-S}.
	The identity~\eqref{eq:equilibrium-RT} follows verbatim from~\cite[Lem.~6]{BCGT:StabilizationfreeHHOPosteriori2023}.
\end{proof}

\begin{theorem}[Reliability and efficiency]
\label{thm:a-posteriori-reliability-efficiency}
There exists $h_0>0$ such that, for $h_{\max}\leq h_0$, the a posteriori error estimator $\eta(\mathcal T)$ from \eqref{eq:local-estimator} satisfies
	\begin{align}
		C_{\mathrm{eff}}^{-2}\eta^2(\mathcal T)
		&\leq
		|\lambda-\lambda_h|
		+\|\mathcal Iu-u_h\|_{a_h}^2
		+\|\nabla u-p_h\|_{L^2}^2
		+\mathrm{osc}_{k+1}^2(Vu,\mathcal T)
		\notag\\
		&\leq
		C_{\mathrm{rel}}^2\eta^2(\mathcal T).
		\label{eq:a-posteriori-efficiency}
	\end{align}
	The constants $h_0,C_{\mathrm{rel}},C_{\mathrm{eff}}>0$ depend only on $\Omega$, the shape regularity, $k$, $\sigma$, $\alpha$, $V$, the constants in \Cref{thm:a_priori}, the eigenvalue $\lambda$, and its spectral separation.
\end{theorem}

\begin{proof}[Proof of \cref{thm:a-posteriori-reliability-efficiency}]
	We prove the upper and lower bounds in~\eqref{eq:a-posteriori-efficiency} separately.

	\medskip\noindent\emph{Step 1: reliability.} Recall from~\Cref{thm:a_priori} and~\eqref{eq:Err} for simple $\lambda$ that
	\begin{align}\label{eqn:a_priori_single}
	\begin{split}
		&|\lambda-\lambda_h|
		+\|\mathcal Iu-u_h\|_{a_h}^2
		+ h_{\max}^{-2s}\|u-u_{\mathcal T}\|_{L^2}^2 + h_{\max}^{-2s_V}\||V|^{1/2}(u-u_{\mathcal T})\|_{L^2}^2\\
		&\qquad\lesssim \operatorname{Err}(u),
	\end{split}
	\end{align}
	provided $h_{\max}\leq h_0$ is sufficiently small. Since $\|u-u_{\mathcal T}\|_{L^2}^2=2-2(u,u_{\mathcal T})_{L^2}$, the unit eigenvector from \cref{thm:a_priori} satisfies $(u,u_{\mathcal T})_{L^2}\geq0$ after possibly reducing $h_0$ and thus coincides with $u$.
	As the elliptic projection satisfies $\|\nabla_{\mathrm{pw}}(u - \mathcal{G} u)\|_{L^2}\leq \|\nabla u-p_h\|_{L^2}$, it suffices to prove reliability of the form
	\begin{align}\label{eqn:short_rel}
		\|\nabla u-p_h\|_{L^2}^2 +\mathrm{osc}_{k+1}^2(Vu,\mathcal T) \lesssim\eta^2(\mathcal T).
	\end{align}

	This proof of~\eqref{eqn:short_rel} follows that of~\cite[Thm.~6.2]{CGT:AdaptiveHybridHighorder2024} and focuses on the potential and correction terms.
	The solution $w\in H_0^1(\Omega)$ to
	$(\nabla w,\nabla v)_{L^2}=(p_h,\nabla v)_{L^2}$ for all $v\in H_0^1(\Omega)$
	leads to the orthogonal decomposition of the gradient error as
	\begin{align}
		\|\nabla u-p_h\|_{L^2}^2
		=
		\|\nabla(u-w)\|_{L^2}^2
		+\|\nabla w-p_h\|_{L^2}^2.
		\label{eq:a-posteriori-helmholtz-split}
	\end{align}
	The final term in~\eqref{eq:a-posteriori-helmholtz-split} is the nonconformity error. As in~\cite[Lem.~5]{BCGT:StabilizationfreeHHOPosteriori2023}, the orthogonality to the divergence-free lowest-order Raviart--Thomas functions from~\eqref{eq:equilibrium-RT} in \Cref{lem:a-posteriori-equilibrium} and a Helmholtz decomposition control it by the tangential jumps in $\eta(\mathcal T)$, i.e.,
	\begin{align}\label{eqn:nonconforming_error}
		\|\nabla w-p_h\|_{L^2}=\min_{v\in H^1_0(\Omega)}\|\nabla v- p_h\|_{L^2}\lesssim 
	\eta(\mathcal T).
	\end{align}
	Note that the piecewise curl contribution from the general estimate in~\cite{BCGT:StabilizationfreeHHOPosteriori2023} vanishes here because $p_h=\nabla_\pw \mathcal{R}u_h$ is an elementwise gradient.

	To control the conforming error $\varrho\coloneqq u-w\in H^1_0(\Omega)$ in~\eqref{eq:a-posteriori-helmholtz-split},
	let $\varrho_C\in S_0^1(\mathcal T)$ denote its Scott--Zhang interpolant.
	Straightforward manipulations for the continuous eigenvalue problem~\eqref{eq:weak_form} and the residual equation~\eqref{eq:equilibrium-S} in \Cref{lem:a-posteriori-equilibrium} imply
	\begin{align}
        \begin{split}
		\|\nabla\varrho\|_{L^2}^2
		={}&
		\big((\lambda u-\lambda_hu_{\mathcal T},\varrho)_{L^2}
		-(V(u-u_{\mathcal T}),\varrho)_{L^2}\big)
		\\
		&+
		\big(((\lambda_h-V)u_{\mathcal T},\varrho-\varrho_C)_{L^2}
		-(p_h,\nabla(\varrho-\varrho_C))_{L^2}\big)
		\\
		&-
		\alpha
		(h_{\mathcal{T}} Q(Vu_{\mathcal T}),h_{\mathcal{T}} Q(V\varrho_C))_{L^2}.
                \end{split}
		\label{eq:reliability-identity}
	\end{align}
	For the first term on the right-hand side of~\eqref{eq:reliability-identity}, a Friedrichs inequality and the multiplier bound~\eqref{eqn:multiplier_12} for $|V|^{1/2}$ provide the upper bound
	\[
	\bigl(
	\|\lambda u-\lambda_hu_{\mathcal T}\|_{L^2}
	+\|V\|_{L^p}^{1/2}\||V|^{1/2}(u-u_{\mathcal T})\|_{L^2}
	\bigr)\|\nabla\varrho\|_{L^2}.
	\]
	Piecewise integration by parts and the Scott--Zhang estimates control the second term in~\eqref{eq:reliability-identity} by $\eta(\mathcal T)\|\nabla\varrho\|_{L^2}$.
	Since $Q$ is $L^2$-orthogonal, the multiplier bound~\eqref{eq:continuous-multiplier} and the $H^1$-stability of the Scott--Zhang interpolation imply
	\begin{align*}
		\|Q(V\varrho_C)\|_{L^2}\lesssim \|\nabla\varrho\|_{L^2}.
	\end{align*}
	Hence the third term in~\eqref{eq:reliability-identity} contributes a small oscillation term
	\begin{align*}
	-\alpha
		(h_{\mathcal{T}} Q(Vu_{\mathcal T}),h_{\mathcal{T}} Q(V\varrho_C))_{L^2}
		\lesssim \alpha h_{\max}\,\mathrm{osc}_{k+1}(Vu_{\mathcal T},\mathcal T)\|\nabla\varrho\|_{L^2}.
	\end{align*}
	The previously discussed estimates for the right-hand side of~\eqref{eq:reliability-identity} establish
	\begin{align}\label{eq:preliminary-reliability}
	\begin{split}
		\|\nabla \varrho\|_{L^2}^2
		&\lesssim
		\eta^2(\mathcal T)
		+\|\lambda u-\lambda_hu_{\mathcal T}\|_{L^2}^2\\
		&\quad+\||V|^{1/2}(u-u_{\mathcal T})\|_{L^2}^2
		+h_{\max}^2\,\mathrm{osc}_{k+1}^2(Vu_{\mathcal T},\mathcal T),
	\end{split}
	\end{align}
	with a hidden constant that depends on $\alpha$ and $V$.

	Standard algebraic manipulations with $\|u\|_{L^2}=1=\|u_{\mathcal T}\|_{L^2}$ reveal
	\begin{align*}
		\|\lambda u-\lambda_hu_{\mathcal T}\|_{L^2}^2
		&=(\lambda-\lambda_h)^2+\lambda\lambda_h\|u-u_{\mathcal T}\|_{L^2}^2.
	\end{align*}
	Hence~\eqref{eqn:a_priori_single} with $\operatorname{Err}(u)\lesssim h_{\max}^{2s}$ from~\eqref{eqn:Err_bound} with $0< s_V\leq s$ leads to
	\begin{align}\label{eqn:a_posteriori_high_order}
		\|\lambda u-\lambda_hu_{\mathcal T}\|_{L^2}^2
		+\||V|^{1/2}(u-u_{\mathcal T})\|_{L^2}^2
		\lesssim h_{\max}^{2s_V}\operatorname{Err}(u).
	\end{align}
	This,~\eqref{eq:a-posteriori-helmholtz-split}--\eqref{eqn:nonconforming_error}, and \eqref{eq:preliminary-reliability} provide
	\begin{align}\label{eqn:reliability_wo_oscillation}
		\|\nabla u-p_h\|_{L^2}^2
			 &\lesssim \eta^2(\mathcal T)
			 + h_{\max}^{2s_V} \operatorname{Err}(u)
		+h_{\max}^2\,\mathrm{osc}_{k+1}^2(Vu_{\mathcal T},\mathcal T).
	\end{align}
	It remains to control the potential oscillation $\mathrm{osc}_{k+1}(Vu,\mathcal T)$.
	The (reverse) triangle inequalities, the $L^2$-stability of $Q$, and $u-u_{\mathcal T}=Qu+(\mathcal Iu-u_h)_{\mathcal T}$ imply
	\begin{align}
		|
		\mathrm{osc}_{k+1}(Vu,\mathcal T)
		&-\mathrm{osc}_{k+1}(Vu_{\mathcal T},\mathcal T)
		|
		\leq
		\mathrm{osc}_{k+1}(V(u-u_{\mathcal T}),\mathcal T)
		\nonumber\\
		&\leq
		\|h_{\mathcal T}VQu\|_{L^2}
		+\|h_{\mathcal T}Q(V(\mathcal Iu-u_h)_{\mathcal T})\|_{L^2}
		\nonumber\\
		&\lesssim
		h_{\max}^{2-n/p}\|\nabla_{\mathrm{pw}}(u-\mathcal Gu)\|_{L^2}
		+h_{\max}^{\min\{1,\,2-n/p\}}\|\mathcal Iu-u_h\|_{a_h},\label{eqn:a_posteriori_V_bound}
	\end{align}
	with \Cref{lem:higher-order-defect,lem:local_pro_err} in the last step.
	This and~\eqref{eqn:a_priori_single} lead to
	\begin{align}\label{eqn:oscillation_swap}
		\mathrm{osc}_{k+1}^2(Vu,\mathcal T)
		&\lesssim\mathrm{osc}_{k+1}^2(Vu_{\mathcal T},\mathcal T) + h_{\max}^{2\min\{1,\,2-n/p\}}\operatorname{Err}(u).
	\end{align}
	Recall $\|\nabla_{\mathrm{pw}}(u - \mathcal{G} u)\|_{L^2}\leq \|\nabla u-p_h\|_{L^2}$, and observe from $\operatorname{div}_\pw p_h+\lambda_hu_{\mathcal T}\in\mathbb P^{k+1}(\mathcal T)$ and the local mesh size equivalence $|T|^{1/n}\approx h_T$ that the volume residual in $\eta(\mathcal T)$ controls $\mathrm{osc}_{k+1}(Vu_{\mathcal T},\mathcal T)\lesssim\eta(\mathcal T)$.
	Thus~\eqref{eqn:reliability_wo_oscillation} and~\eqref{eqn:oscillation_swap} with $s_V\leq 1-n/(2p)\leq\min\{1,2-n/p\}$ lead to
	\begin{align*}
		\|\nabla u-p_h\|_{L^2}^2+\mathrm{osc}_{k+1}^2(Vu,\mathcal T)
			 &\lesssim \eta^2(\mathcal T)
		+ h_{\max}^{2s_V} \left(\|\nabla u-p_h\|_{L^2}^2+\mathrm{osc}_{k+1}^2(Vu,\mathcal T)\right).
	\end{align*}
	By possibly reducing $h_0$ and thereby $h_{\max}\leq h_0$, the final sum can be absorbed into the left-hand side to prove~\eqref{eqn:short_rel}, which concludes the proof of reliability.

	\medskip\noindent\emph{Step 2: efficiency.} Standard bubble function arguments (as in~\cite{Ver96a} and specifically for HHO in~\cite{BCGT:StabilizationfreeHHOPosteriori2023,CGT:AdaptiveHybridHighorder2024}) applied with the exact source $f\coloneqq (\lambda-V)u$ and the polynomial discrete source $f_h\coloneqq\lambda_hu_{\mathcal T}-\Pi_{\mathcal T}^{k+1}(Vu_{\mathcal T})$ provide
\begin{align}
	\eta(\mathcal T)
	\lesssim{}&
	\|\nabla u-p_h\|_{L^2}
	+\|h_{\mathcal T}(\lambda u-\lambda_hu_{\mathcal T})\|_{L^2}
	\nonumber\\
	&+\|h_{\mathcal T}V(u-u_{\mathcal T})\|_{L^2}
	+\mathrm{osc}_{k+1}(Vu_{\mathcal{T}},\mathcal T).
	\label{eqn:preliminary_efficiency}
\end{align}
The last three terms in~\eqref{eqn:preliminary_efficiency} originate from the source mismatch $f-f_h=(\lambda u-\lambda_hu_{\mathcal T})-V(u-u_{\mathcal T})-Q(Vu_{\mathcal T})$.
The positive powers of \(h_{\max}\leq h_0\) are collected in the hidden constant so that~\eqref{eqn:a_posteriori_high_order}, \eqref{eqn:a_posteriori_V_bound}, and \cref{lem:higher-order-defect,lem:local_pro_err} show
\begin{align*}
	\|h_{\mathcal T}(\lambda u-\lambda_hu_{\mathcal T})\|_{L^2}
	&+\|h_{\mathcal T}V(u-u_{\mathcal T})\|_{L^2}
	+\mathrm{osc}_{k+1}(Vu_{\mathcal{T}},\mathcal T)\\
	&\lesssim
	\|\nabla u-p_h\|_{L^2}
	+\|\mathcal Iu-u_h\|_{a_h}
	+\mathrm{osc}_{k+1}(Vu,\mathcal T).
\end{align*}
This and~\eqref{eqn:preliminary_efficiency} conclude the proof of efficiency.
\end{proof}

\section{Numerical experiments}%
\label{sec:Numerical experiments}
\begin{figure}[]
	\noindent\hspace{-.2em}\fbox{\parbox{.985\textwidth}{
		\noindent\textbf{Adaptive algorithm (AFEM).}\\[.2em]
		\Input Initial triangulation $\T_0$ and index $m\in\mathbb N$ of the eigenvalue
		$\lambda$
\begin{description}[leftmargin=\widthof{\Keyword{for}\,},topsep=0pt]
	\item[\Keyword{for}] $\ell=0,1,2,\dots$
	\begin{description}[leftmargin=1em]
		\item[\Keyword{SOLVE}] \eqref{eq:discrete_eigenproblem_modif} for
		$m$-th eigenpair $(\lambda_\ell,u_\ell)$
			on $\T_\ell$
			and compute $\gamma_\ell=\gamma(\lambda_\ell)$ from \eqref{eqn:gamma_def}
		\item[\Keyword{ESTIMATE}] local error indicators $\eta_\ell(T)$
			from~\eqref{eq:local-estimator} for all $T\in\mathcal T_\ell$
		\begin{description}[leftmargin=1em]
			\item[\Keyword{if}] $\gamma_\ell=0$ \Keyword{and} $\lambda_\ell\geq0$ \Keyword{then} $\mathrm{GLB}_\ell=\lambda_\ell\leq \lambda$
				is a guaranteed lower bound
			\begin{description}[leftmargin=1em]
				\item[\Keyword{MARK}] minimal subset $\mathcal M_\ell\subset\T_\ell$ with
					$\frac{1}{2}\,\eta_\ell^2(\T_\ell) \leq \sum^{}_{T\in\mathcal{M}_\ell} \eta_\ell^2(T)$
				\item[\Keyword{REFINE}] $\T_\ell$ to smallest NVB refinement $\T_{\ell+1}$ with
					$\mathcal M_\ell\subseteq\T_\ell\setminus\T_{\ell+1}$
			\end{description}
			\item[\Keyword{else}] $\mathrm{GLB}_\ell=c_A(\lambda_\ell)_+/(c_A+\gamma_\ell)\leq \lambda$
				is a guaranteed lower bound
			\begin{description}[leftmargin=1em]
				\item[\Keyword{REFINE}] $\T_\ell$ by uniform red refinement to $\T_{\ell+1}$
			\end{description}
		\end{description}
	\end{description}
\end{description}
\Keyword{Output:} Sequences of triangulations $\T_\ell$, error estimators 
$\eta_\ell=\eta_\ell(\mathcal T_\ell)$,
guaranteed lower bounds $\mathrm{GLB}_\ell\leq \lambda$, and discrete solutions $u_\ell$
}}
	\caption{Mesh refinement strategy adapted from~\cite{CGP:ComparisonGuaranteedLower2025}.}
	\label{fig:AFEM}
\end{figure}

Three computational benchmarks for piecewise-constant, polynomial, and unbounded potentials validate the GLB property and provide striking evidence for higher-order convergence rates.

\begin{figure}
\begin{center}
	\includegraphics[width=.3\textwidth]{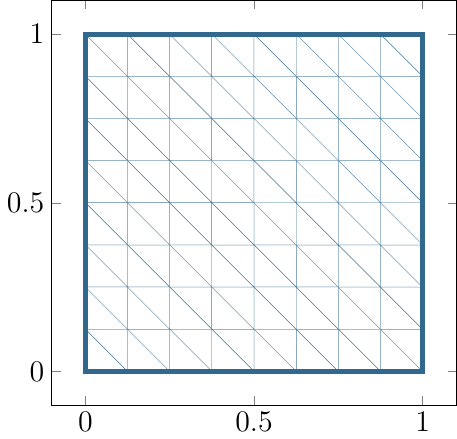}
	\includegraphics[width=.3\textwidth]{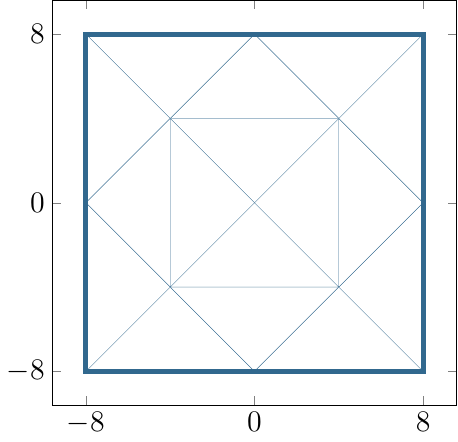}
	\includegraphics[width=.3\textwidth]{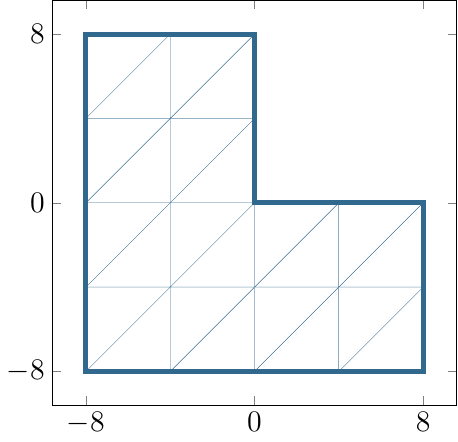}
\end{center}
\caption{Initial triangulations of the unit square (left), the big square $(-8,8)^2$ (middle), and the big L-shaped domain $(-8,8)^2\setminus[0,8]^2$ (right) from~\cite{Carstensen2026}.}
\label{fig:initial_mesh}
\end{figure}

\subsection{Computational setup and parameter choices}%
\label{sub:Computational setup and parameter choices} The MATLAB implementation is based on the code from~\cite{CGT:AdaptiveHybridHighorder2024,Tran2024} and solves the algebraic eigenvalue problem with \texttt{eigs} up to machine precision.
We employ quadrature rules that integrate polynomials exactly. 
On elements containing a logarithmic singularity in \Cref{sub:Logarithmic}, we combine Gaussian quadrature with a power-graded Duffy transformation~\cite{Duffy1982}.
Numerical errors from inexact solves and quadrature are neglected.

Adaptive mesh refinement is compared with uniform red refinement.
\Cref{fig:AFEM} displays the natural adaptive algorithm for these skeletal GLB, based on newest vertex bisection (NVB) and Dörfler marking driven by the error estimators~\eqref{eq:local-estimator}.
All triangulations are refinements of the initial ones from \cref{fig:initial_mesh} and consist exclusively of right-isosceles triangles.
The convergence history plots display the GLB error $\lambda-\mathrm{GLB}$ or the squared error estimator $\eta^2(\mathcal{T})$ for the corresponding eigenvalue index with respect to the number $|\mathcal{T}|$ of mesh elements.
All displayed computations satisfy the lower bound $\mathrm{GLB}\leq\lambda$ for the stated reference value of $\lambda$.

Numerical approximations of the optimal $\Lambda, C_P$, and $C_{\mathrm{st}}$, computed as in~\cite[Sec.~4.2]{CGZ:HigherorderGuaranteedLower2026}, are given in \cref{tab:Lambda} and enter the canonical choice $\sigma = \Lambda^{-1}/3$ and $\alpha = 3 C_P^2$ from \cref{rem:parameters}.
Nonnegative potentials in \Cref{sub:Disorder,sub:Harmonic} lead to the natural choices $c_A=1$ and $\delta_h=0$.
\Cref{sub:Logarithmic} calculates $c_A$ and $\delta_h$ for a nonzero negative potential part $V_-\in L^3(\Omega)$ with the auxiliary constants for~$p=3=r$,
\begin{align*}%
	C_{S,3}=\frac{|\Omega|^{1/3}3^{5/4}}{2^{13/6}\pi}
	\qquad\text{and}\qquad
C_{P,3}^2 = \frac{3}{2}\left(\frac{2}{\pi}\right)^{1/3} C_P^{4/3} C_{\mathrm{st}}^{2/3},
\end{align*}
from~\cite[Cor.~A.2]{TSMO:SharpNumericalInclusion2017} and from \Cref{lem:Lr_right_iso} with $C_P,C_{\mathrm{st}}$ from~\Cref{tab:Lambda}.%
\begin{table}[t]  
\centering  
\begin{tabular}{cccc}
	\hline
	$k$ & $\Lambda$ & $C_{P}^2$ & $C_{\mathrm{st}}^2$
	\\ \hline
	0     & 1.5765e-02 & 1.5416e-02 & 1.5428\\
	1     & 8.8171e-03 & 7.4771e-03 & 1.7280\\
	2     & 6.7597e-03 & 4.6327e-03 & 2.0729\\
	3     & 4.8495e-03 & 3.2398e-03 & 2.3987\\
	4     & 3.7720e-03 & 2.4357e-03 & 2.7506\\
	5     & 2.9469e-03 & 1.9170e-03 & 3.1290
	\end{tabular}
	\caption{Values of $\Lambda$, $C_{P}^2$, and $C_{\mathrm{st}}^2$ for right-isosceles triangles.}
\label{tab:Lambda}
\end{table}%

\subsection{Disorder potential}\label{sub:Disorder}
The disorder potential $V(x)$ from~\cite[Sec.~7.4]{Carstensen2026} with
\begin{equation*}
	V|_{\Omega_{s,t}} \equiv x_{s,t}\in\{0,\dots,10^4-1\}
\end{equation*}
is piecewise constant over a uniform $8\times 8$ partition of the unit square $\Omega = (0,1)^2$,
\begin{align*}
	\Omega_{s, t} = \left[ \tfrac{s - 1}{8},\, \tfrac{s}{8} \right) \times \left[ \tfrac{t - 1}{8},\, \tfrac{t}{8} \right)
	\qquad\text{for all }s,t =1, \dots, 8,
\end{align*}
with integer coefficient matrix $X\in \mathbb N_0^{8\times8}$ obtained%
\footnote{This benchmark initializes the MATLAB random number generator by \texttt{rng(5, 'twister')}.}
from \texttt{randi(1e4, [8,8])-1}.

Since the initial triangulation displayed in~\cref{fig:initial_mesh} resolves the partition $\{\Omega_{s,t}\}_{s,t}$, the potential is piecewise constant over all computed triangulations.
Consequently, $Q(Vv_{\mathcal T})=0$ for all $v_{\mathcal T}\in\mathbb P^{k+1}(\mathcal T)$ and the negative correction term in~\eqref{eq:modified_bilinear_form} vanishes.
The discrete problem is therefore independent of $\alpha$.
Formally setting $\alpha^{-1}=0$ in \cref{theo:main_lower_bound} recovers the GLB from~\cite{Tran2024}, but the computations here retain the universal parameter choice from~\cref{sub:Computational setup and parameter choices}.

\begin{figure}
	\centering
	\includegraphics[width=.9\textwidth]{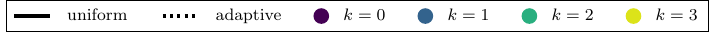}\\
	\hspace*{-3em}
	\hbox{
	\includegraphics[width=0.55\textwidth]{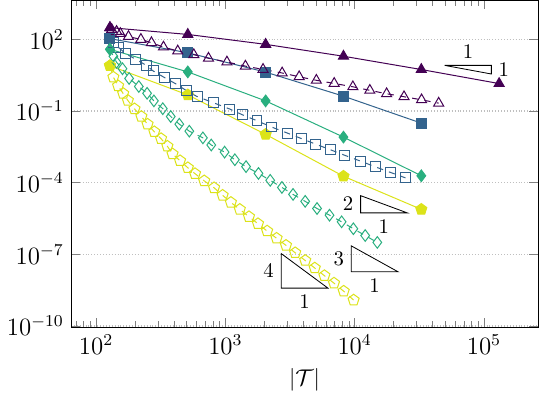}
	\includegraphics[width=0.55\textwidth]{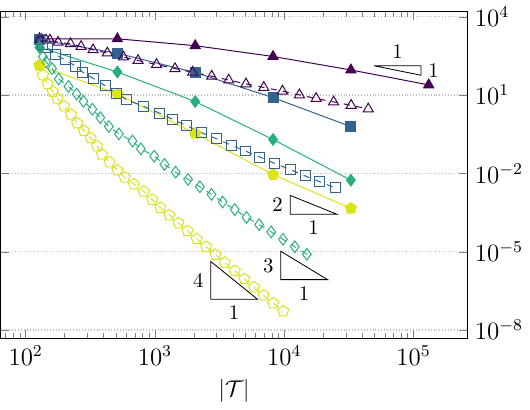}}
	\caption{Convergence history plots of the $\mathrm{GLB}$ error (left) and the error estimator $\eta^2(\mathcal{T})$ (right) for $\lambda_1$ with the disorder potential on the unit square in \Cref{sub:Disorder}.
	}
	\label{fig:Disorder}
\end{figure}

\begin{figure}
	\centering
	\hspace*{-3em}
	\hbox{
	\includegraphics[width=0.45\textwidth]{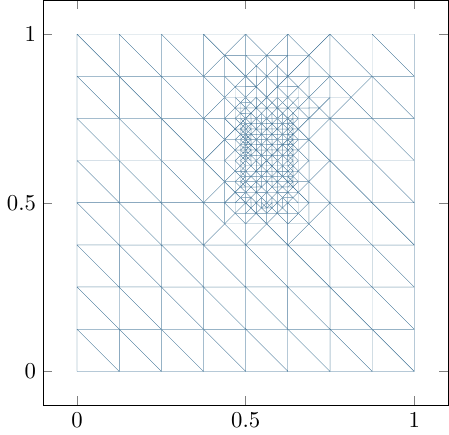}
	\hspace{1em}
	\includegraphics[width=0.45\textwidth]{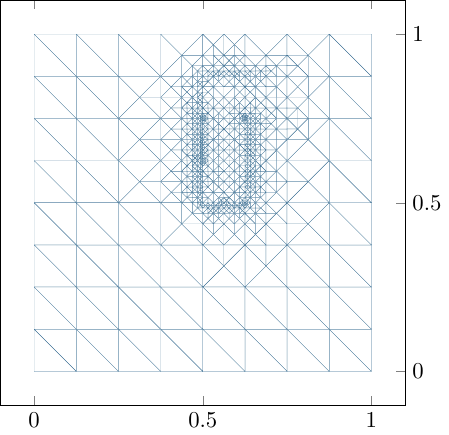}
}
	\caption{Adaptive meshes for $k=0$ (left, $|\T|=740$) and $k=2$ (right, $|\T|=1445$) with the disorder potential on the unit square in \Cref{sub:Disorder}.}
	\label{fig:Anderson_mesh}
\end{figure}

The ground energy $\lambda_{1}=539.16091242$ is obtained from a fine reference computation.
Uniform refinement leads in~\cref{fig:Disorder} to the optimal convergence rate $k+1$ only for $k\leq1$, while adaptive refinement recovers the optimal convergence rates for the GLB error and $\eta^2(\mathcal T)$ for all polynomial degrees $k$ over the displayed range.
The reduced convergence rates under uniform refinement are consistent with the limited piecewise regularity of the eigenfunction at the junctions of the potential jumps.
The strong Anderson localization of the principal eigenfunction explains the smaller errors compared to uniform meshes of comparable cardinality and the local refinement displayed in \cref{fig:Anderson_mesh}, which leaves large parts of the domain coarse.

\subsection{Harmonic potential}\label{sub:Harmonic}
The harmonic potential $V(x) = |x|^2/2$ represents a smooth, convex confinement trap centred at the origin, and is considered on the big square and the big L-shaped domain with the initial triangulations displayed in \Cref{fig:initial_mesh}.
The Gaussian decay of the full-space eigenfunctions explains the weak influence of the Dirichlet condition at the outer boundary $\partial (-8,8)^2$ on the lower eigenvalues.

\begin{figure}
	\centering
	\includegraphics[width=.9\textwidth]{Figures/Legend_k_out.pdf}\\
	\hspace*{-3em}
	\hbox{
	\includegraphics[width=0.55\textwidth]{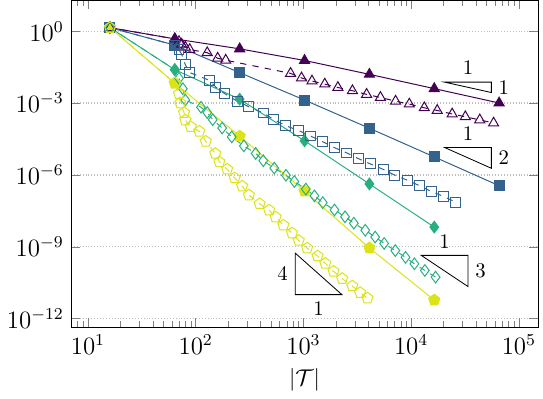}
	\includegraphics[width=0.55\textwidth]{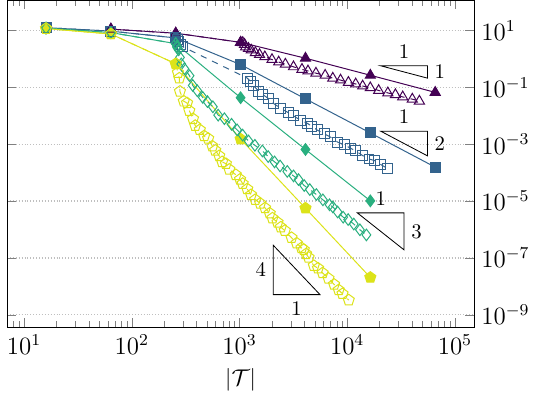}}
	\caption{Convergence history plots of the $\mathrm{GLB}$ error for $\lambda_1$ (left) and $\lambda_{50}$ (right) with the harmonic potential on the big square in \Cref{sub:Harmonic}.
	}
	\label{fig:HarmonicSquare}
\end{figure}

The full-space eigenvalues $j\sqrt2$ have multiplicity $j$ for $j\in\mathbb N$.
Zero extension, the min--max principle, and unique continuation imply that the Dirichlet eigenvalues on $\Omega$ are strictly larger than their full-space counterparts.
On the big square, the ground energy agrees with its full-space counterpart $\sqrt2$ to machine precision, while $\lambda_{50}=14.14213562490$ belongs to a cluster of five distinct double eigenvalues that lie within $(10\sqrt2,10\sqrt2+10^{-7})$ around the ten-fold full-space eigenvalue $10\sqrt2$.
This cluster has been confirmed by disjoint enclosures obtained from the high-order GLB and conforming upper bounds.
\Cref{theo:main_lower_bound} guarantees the GLB on every mesh, so only the refinement for the double eigenvalue $\lambda_{50}$ goes beyond \Cref{sec:a-posteriori}, which considers simple eigenvalues.
\Cref{fig:HarmonicSquare} displays the optimal convergence order $k+1$ of the GLB for $\lambda_1$ and $\lambda_{50}$ under uniform and adaptive refinement.
Adaptivity improves the error substantially for larger $k$, already on coarse meshes.

\begin{figure}
	\centering
	\includegraphics[width=.9\textwidth]{Figures/Legend_k_out.pdf}\\
	\hspace*{-3em}
	\hbox{
	\includegraphics[width=0.55\textwidth]{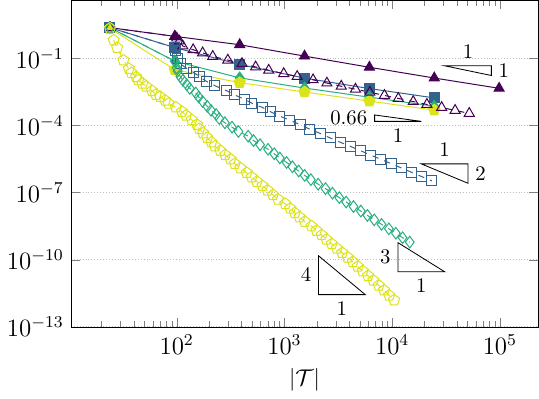}
	\includegraphics[width=0.55\textwidth]{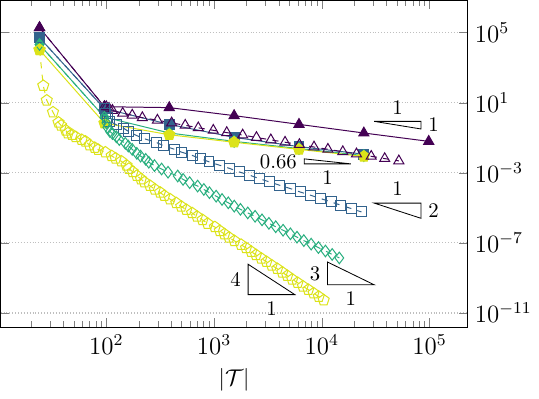}}
	\caption{Convergence history plots of the $\mathrm{GLB}$ error (left) and $\eta^2(\mathcal{T})$ (right) for $\lambda_1$ with the harmonic potential on the L-shaped domain in \Cref{sub:Harmonic}.
	}
	\label{fig:HarmonicLshape}
\end{figure}

\begin{figure}
	\centering
	\hspace*{-3em}
	\hbox{
	\includegraphics[width=0.45\textwidth]{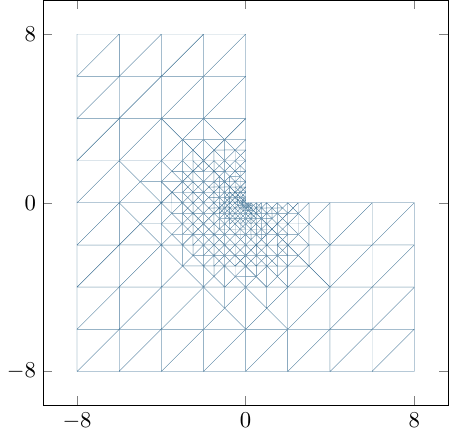}
	\hspace{1em}
	\includegraphics[width=0.45\textwidth]{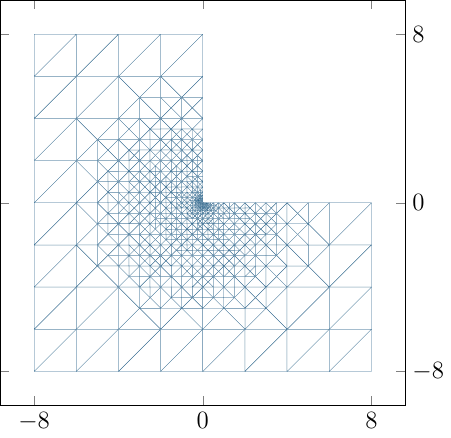}
}
	\caption{Adaptive meshes for $k=0$ (left, $|\T|=714$) and $k=2$ (right, $|\T|=1482$) for the harmonic potential on the big L-shape in \Cref{sub:Harmonic}.}
	\label{fig:Harmonic_mesh}
\end{figure}

On the big L-shaped domain, the re-entrant corner at the origin results in a reduced observed convergence rate $0.66$ of the GLB and the error estimator $\eta^2(\mathcal{T})$ for the principal eigenvalue $\lambda_1=2.357022603953$ displayed in \Cref{fig:HarmonicLshape}.
Adaptive refinement as shown in \cref{fig:Harmonic_mesh} recovers the optimal convergence rate.

All displayed examples for both domains produce a negative $\lambda_{h}<0$ only on the coarsest mesh, except for $k=3$ on the big L-shaped domain, where $0<\operatorname{GLB}=\lambda_h\leq \lambda$ holds on all meshes.
Possible further uniform refinements in the AFEM algorithm of \Cref{fig:AFEM} are then triggered by $\gamma(\lambda_{h})>0$.


\subsection{Logarithmic potential}%
\label{sub:Logarithmic}
Consider the unbounded logarithmic potential
\begin{align*}
	V(x) = \frac{3}{64}\log\left(\frac{8\sqrt2}{|x-z_+|}\right)-\frac{1}{8}\log\left(\frac{1}{|x-z_-|}\right)
\end{align*}
with poles at $z_+=(0,0)$ and $z_-=(-4,-4)$ on the big square domain, satisfying $V\in L^p(\Omega)$ for any $p<\infty$.
Its negative part is supported in the unit ball about $z_-$ with majorant $\tfrac{1}{8}\log(|x-z_-|^{-1})$.
For $p=3$, integration in polar coordinates gives
\begin{align*}
	\|V_-\|_{L^3}\leq \frac{\|\log(|x|)\|_{L^3(|x|<1)}}{8}=\frac{(3\pi/4)^{1/3}}{8}\eqqcolon C_{-}.
\end{align*}
Hence \Cref{lem:Garding_Lp,lem:negative_potential_Lp} lead with the values $C_{S,3}$ and $C_{P,3}$ from~\Cref{sub:Computational setup and parameter choices}~to
\begin{align*}
	\delta_h^2=\frac{3^{4/3}}{2^{13/3}}C_P^{4/3}C_{\mathrm{st}}^{2/3} h_{\max}^{4/3}
	\qquad\text{and}\qquad
	c_A=1-\frac{3^{17/6}}{2^{8/3}\pi^{5/3}}=0.47457>0.
\end{align*}

\begin{figure}[]
    \centering
	\includegraphics[width=.9\textwidth]{Figures/Legend_k_out.pdf}\\
    \hspace*{-3em}
    \hbox{
    \includegraphics[width=0.55\textwidth]{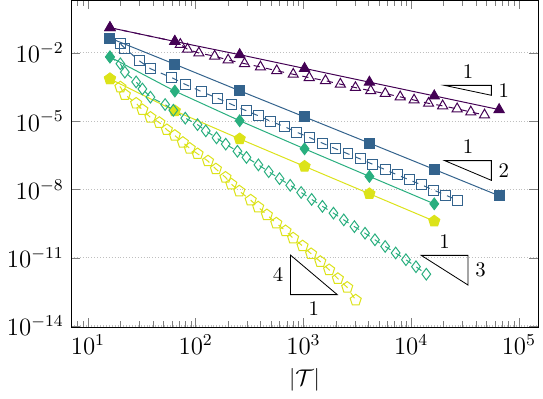}
    \includegraphics[width=0.55\textwidth]{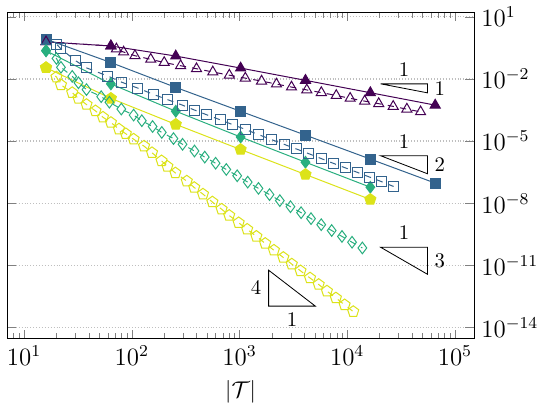}
    }
	\caption{Convergence histories of the GLB error (left) and $\eta^2(\mathcal{T})$ (right) for $\lambda_1$ with the logarithmic
        potential in \Cref{sub:Logarithmic}.}
    \label{fig:Logarithmic}
\end{figure}

\begin{figure}
    \centering
    \hspace*{-3em}
    \hbox{
    \includegraphics[width=0.45\textwidth]{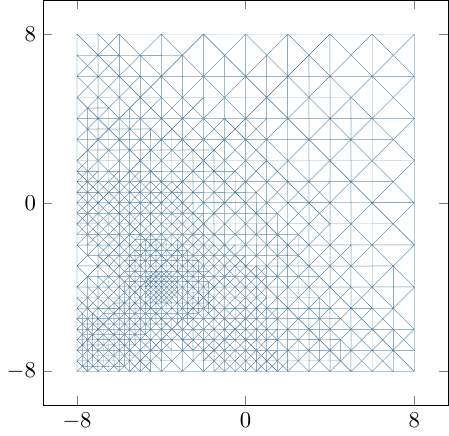}
	\hspace{1em}
    \includegraphics[width=0.45\textwidth]{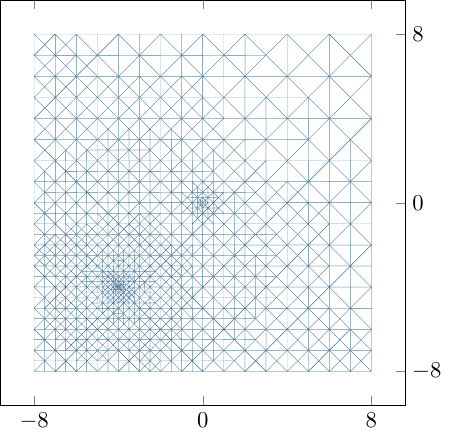}}
    \caption{Adaptive meshes for $k=0$ (left, $|\T|=2202$) and $k=2$
        (right, $|\T|=1900$) with the logarithmic potential in \Cref{sub:Logarithmic}.}
    \label{fig:LogarithmicSquare_mesh}
\end{figure}

\Cref{fig:Logarithmic} displays reduced convergence rates under uniform refinement for $k\geq2$ towards the ground energy $\lambda_1=0.299567750934152$, whereas adaptive refinement recovers the optimal rates $k+1$ for the GLB error and the error estimator $\eta^2(\mathcal T)$.
The singularities cause no visible reduction in the uniform convergence rate for $k=0,1$ as expected from $u\in H^{3-\varepsilon}(\Omega)$ for any eigenfunction and $\varepsilon>0$ by elliptic regularity.
The adaptive mesh for $k=0$ in \cref{fig:LogarithmicSquare_mesh} is only refined towards the attractive point $z_-$ where the mass of the eigenfunction concentrates while the repulsive singularity at the origin $z_+$ does not trigger local refinement.
For larger polynomial degree $k\geq2$, adaptive refinement towards both singularities $z_{\pm}$ recovers optimal convergence rates.
The direct GLB condition $\gamma(\lambda_{h,1})=0$ holds on all displayed meshes except for $k=0$ on the initial triangulation, which provides the lower bound $\operatorname{GLB}\leq \lambda$ via postprocessing by $\gamma(\lambda_h)>0$.

\FloatBarrier
\paragraph{Acknowledgements}
The authors would like to thank T.~Tran (University of Augsburg) for providing an implementation of the HHO method on which our code is based. M.~Hauck acknowledges funding from the Deutsche Forschungsgemeinschaft (DFG, German Research Foundation) -- Project-ID 258734477 -- SFB 1173.

\paragraph{Use of large language models}
ChatGPT and Claude were used for language editing, coding assistance, literature review, and to check proof arguments. The authors have checked all mathematical arguments and take full responsibility for the content of the paper.

\bibliography{bib,Numerik}
\bibliographystyle{alphaabbr}

\appendix
\FloatBarrier
\section{A priori error estimate of the source problem}\label{sec:apriori:A}

\begin{lemma}[Source problem error estimate]\label{lem:err_est:source}
For sufficiently small $h_{\max}$ and any $f\in L^2(\Omega)$, let $u=Tf$ denote the solution of~\cref{eq:source_pro} and $u_h=S_hf = (u_{\T},u_{\F})$ the corresponding  discrete approximation defined by~\cref{eq:dis_source_pro}. Then, the following estimates hold:
	\begin{align*}
		&\|\mathcal{I} u - u_h\|_{a_h} + \||V|^{1/2}(u_{\mathcal{T}} - \Pi_{\mathcal{T}}^{k+1}u)\|_{L^2} + |u_h|_{s_h} \\
		&\qquad \lesssim \|\nabla_{\mathrm{pw}}(u - \mathcal{G} u)\|_{L^2} + \mathrm{osc}_{k+1}(f, \mathcal{T}) + \mathrm{osc}_{k+1}(Vu, \mathcal{T}).
	\end{align*}
\end{lemma}

\begin{proof}
Let $e_h\coloneqq\mathcal{I} u-u_h$. Testing the weak formulation~\cref{eq:source_pro} with $\mathcal{J}e_h$ and the discrete formulation~\cref{eq:dis_source_pro} with $e_h$, and subtracting, yields
	\begin{equation}\label{eq:priori_source_auxi1}
		\begin{aligned}
			&(\nabla u, \nabla \mathcal{J}e_h)_{L^2} - (\nabla_{\mathrm{pw}}\mathcal{R}u_h, \nabla_{\mathrm{pw}}\mathcal{R}e_h)_{L^2} + (Vu, \mathcal{J}e_h)_{L^2} - (Vu_{\mathcal{T}}, e_{\mathcal{T}})_{L^2} \\
			&\qquad = (f, \mathcal{J}e_h - e_{\mathcal{T}})_{L^2} + s_h(u_h, e_h) - \alpha (h_{\mathcal{T}} Q(Vu_{\mathcal{T}}), h_{\mathcal{T}} Q(Ve_{\mathcal{T}}))_{L^2}.
		\end{aligned}
	\end{equation}
	To rewrite the reconstructed gradient contribution, we use the property from~\cref{eq:keyidentity} and~\cref{eq:defGh}
	$$(\nabla u, \nabla_{\mathrm{pw}}\mathcal{R}e_h)_{L^2} = (\nabla_{\mathrm{pw}}\mathcal{R}\mathcal{I} u, \nabla_{\mathrm{pw}}\mathcal{R}e_h)_{L^2}$$
	 to expand this term as
	\begin{align*}
		&-(\nabla_{\mathrm{pw}}\mathcal{R}u_h, \nabla_{\mathrm{pw}}\mathcal{R}e_h)_{L^2}\\
		 &\qquad = -(\nabla u, \nabla_{\mathrm{pw}}\mathcal{R}e_h)_{L^2} + \|\nabla_{\mathrm{pw}}\mathcal{R}e_h\|_{L^2}^2 \\
		&\qquad = -(\nabla u, \nabla_{\mathrm{pw}}(\mathcal{R}e_h - \mathcal{J}e_h))_{L^2} - (\nabla u, \nabla \mathcal{J}e_h)_{L^2} + \|\nabla_{\mathrm{pw}}\mathcal{R}e_h\|_{L^2}^2.
	\end{align*}
	Similarly, we decompose the potential term, the stabilization term, and the correction term in \cref{eq:priori_source_auxi1} as follows:
	\begin{align*}
		(Vu, \mathcal{J}e_h)_{L^2} - (Vu_{\mathcal{T}}, e_{\mathcal{T}})_{L^2} &= (Vu, \mathcal{J}e_h - e_{\mathcal{T}})_{L^2} + (V e_{\mathcal{T}}, e_{\mathcal{T}})_{L^2}\\
		&\quad  + (VQ u, e_{\mathcal{T}})_{L^2}, \\
		s_h(u_h, e_h) &= \tfrac{1}{2}\left(|\mathcal{I} u|_{s_h}^2 - |e_h|_{s_h}^2 - |u_h|_{s_h}^2\right), \\
		-\alpha(h_{\mathcal{T}}Q(Vu_{\mathcal{T}}), h_{\mathcal{T}}Q(Ve_{\mathcal{T}}))_{L^2} &= -\alpha(h_{\mathcal{T}}Q(Vu), h_{\mathcal{T}}Q(Ve_{\mathcal{T}}))_{L^2} \\
		&\quad + \alpha(h_{\mathcal{T}}Q(VQ u), h_{\mathcal{T}}Q(Ve_{\mathcal{T}}))_{L^2} \\
		&\quad + \alpha(h_{\mathcal{T}}Q(Ve_{\mathcal{T}}), h_{\mathcal{T}}Q(Ve_{\mathcal{T}}))_{L^2}.
	\end{align*}
	Substituting these expansions back into \cref{eq:priori_source_auxi1} and moving the norms of the error to the left-hand side leads to the following identity:
	\begin{align}
		\|\nabla_{\mathrm{pw}}\mathcal{R}e_h\|_{L^2}^2 &+ (V e_{\mathcal{T}},e_{\mathcal{T}})_{L^2} + \tfrac{1}{2}\left(|e_h|_{s_h}^2 + |u_h|_{s_h}^2\right) - \alpha\|h_{\mathcal{T}} Q(Ve_{\mathcal{T}})\|_{L^2}^2 \notag\\
        \begin{split}\label{eq:rhserrest}
                    ={}& (\nabla u, \nabla_{\mathrm{pw}}(\mathcal{R}e_h - \mathcal{J}e_h))_{L^2} + \tfrac{1}{2}|\mathcal{I} u|_{s_h}^2 + (f - Vu, \mathcal{J}e_h - e_{\mathcal{T}})_{L^2} \\
		&- (VQ u, e_{\mathcal{T}})_{L^2} - \alpha(h_{\mathcal{T}} Q(Vu), h_{\mathcal{T}} Q(Ve_{\mathcal{T}}))_{L^2} \\
		&+ \alpha(h_{\mathcal{T}} Q(VQ u), h_{\mathcal{T}} Q(Ve_{\mathcal{T}}))_{L^2}.
        \end{split}
	\end{align}
We bound each term on the right-hand side of the latter identity individually. For the first term in \cref{eq:rhserrest}, we have
\begin{align*}
    (\nabla u, \nabla_{\mathrm{pw}}(\mathcal{R}e_h - \mathcal{J}e_h))_{L^2}
    &= -(\nabla_{\mathrm{pw}} (u - \mathcal{G} u), \nabla\mathcal{J}e_h)_{L^2} \\
    &\lesssim \|\nabla_{\mathrm{pw}}(u - \mathcal{G} u)\|_{L^2}\|e_h\|_{a_h},
\end{align*}
where the equality utilizes the orthogonality $\nabla_{\mathrm{pw}}(\mathcal{J}e_h - \mathcal{R} e_h) \perp \nabla_{\mathrm{pw}} \mathbb{P}^{k+1}(\mathcal{T})$ from $\mathcal Re_h=\mathcal R\mathcal I\mathcal Je_h=\mathcal G\mathcal Je_h$ and the inequality follows from the stability bound \cref{eq:Jstability} applied to $e_h$.

The second term in \cref{eq:rhserrest} can be estimated as
\begin{equation*}
    |\mathcal{I} u|_{s_h} \lesssim \|\nabla_{\mathrm{pw}}(u - \mathcal{G} u)\|_{L^2},
\end{equation*}
which is the global counterpart of the local bound \cref{eq:proof_stabilization}.

For the third term in \cref{eq:rhserrest}, we obtain
\begin{align*}
    (f - Vu, \mathcal{J}e_h - e_{\mathcal{T}})_{L^2}
    &= (Q(f - Vu), Q\mathcal{J}e_h)_{L^2} \\
    &\lesssim \mathrm{osc}_{k+1}(f - Vu, \mathcal{T}) \|e_h\|_{a_h},
\end{align*}
which relies on the identity $\mathcal{J} e_h - e_{\mathcal{T}} = Q (\mathcal{J} e_h)$, the $L^2$-orthogonality properties of $Q$, \eqref{eqn:CP}, and the stability estimate \cref{eq:Jstability}.

For the remaining terms on the right-hand side of~\cref{eq:rhserrest}, we obtain
\begin{align*}
    (VQ u, e_{\mathcal{T}})_{L^2}
    &= (Q u, Q(Ve_{\mathcal{T}}))_{L^2} \\
    &\lesssim h_{\max}^{\min\{1,\,2-n/p\}}
        \|\nabla_{\mathrm{pw}}(u-\mathcal{G} u)\|_{L^2}
        \|e_h\|_{a_h}, \\[0.5em]
    (h_{\mathcal{T}}Q(Vu),
     h_{\mathcal{T}}Q(Ve_{\mathcal{T}}))_{L^2}
    &\lesssim h_{\max}^{\min\{1,\,2-n/p\}}
        \mathrm{osc}_{k+1}(Vu,\mathcal{T})
        \|e_h\|_{a_h}, \\[0.5em]
    (h_{\mathcal{T}}Q(VQ u),
     h_{\mathcal{T}}Q(Ve_{\mathcal{T}}))_{L^2}
    &\lesssim h_{\max}^{2-n/p+\min\{1,\,2-n/p\}}
        \|\nabla_{\mathrm{pw}}(u-\mathcal{G} u)\|_{L^2}
        \|e_h\|_{a_h}.
\end{align*}
Here, we have used \cref{lem:higher-order-defect,lem:local_pro_err}, as well as the $L^2$-stability of $Q$ and \cref{eqn:CP}.

	Discrete coercivity from \Cref{lem:shifted-coercivity} bounds the left-hand side of~\eqref{eq:rhserrest} from below by $c_0\|e_h\|_{a_h}^2+\frac12|u_h|_{s_h}^2$, while a weighted Young inequality absorbs the contributions involving $\|e_h\|_{a_h}$ from the right-hand side.
The identity $e_{\mathcal T}=\mathcal Je_h-Q\mathcal Je_h$, \Cref{lem:potential_estimates}, and the stability of $\mathcal J$ give $\||V|^{1/2}e_{\mathcal T}\|_{L^2}\lesssim\|e_h\|_{a_h}$.
This concludes the proof.
\end{proof}

To obtain an improved $L^2$-error estimate, we employ an Aubin--Nitsche-type argument, gaining a factor of \(h_{\max}^s\).

\begin{lemma}[Source problem $L^2$-estimate]\label{lem:L2err_est:source}
Consider the setup of \cref{lem:err_est:source} and assume elliptic regularity with index $s \in (0,1]$ as in \cref{eq:ellipticregularity}. Then every $f \in L^2(\Omega)$ satisfies
    \begin{equation*}
        \|\Pi_{\mathcal{T}}^{k+1} u-u_{\mathcal{T}} \|_{L^2} \lesssim h_{\max}^s \left( \|\nabla_{\mathrm{pw}}(u - \mathcal{G} u)\|_{L^2} + \mathrm{osc}_{k+1}(f, \mathcal{T}) + \mathrm{osc}_{k+1}(Vu, \mathcal{T}) \right).
    \end{equation*}
\end{lemma}
\begin{proof}
    Abbreviate $e_h \coloneqq \mathcal{I} u - u_h$ and let $z = T e_{\mathcal{T}} \in H_0^1(\Omega)$ denote the solution to the dual problem \eqref{eq:source_pro}. Using the moment-preserving properties of $\mathcal{J}$ (cf.~\cref{lem:enrichment}) and the definition of the dual problem, we obtain
    \begin{align}
        \|e_{\mathcal{T}}\|_{L^2}^2 &= (e_{\mathcal{T}}, \mathcal{J} e_h)_{L^2}\notag \\
        &= (\nabla z, \nabla \mathcal{J} e_h)_{L^2} + (V z, \mathcal{J} e_h)_{L^2}\notag \\
        &= (\nabla z, \nabla_{\mathrm{pw}}(\mathcal{J} e_h - \mathcal{R} e_h))_{L^2} + (\nabla z, \nabla_{\mathrm{pw}} \mathcal{R} e_h)_{L^2} + (V z, \mathcal{J} e_h)_{L^2}.\label{eq:rhsopL21}
    \end{align}
  For the first term on the right-hand side, the orthogonality relation $\nabla_{\mathrm{pw}}(\mathcal{J}e_h - \mathcal{R} e_h) \perp \nabla_{\mathrm{pw}} \mathbb{P}^{k+1}(\mathcal{T})$, the definition \cref{eq:defGh} of the elliptic projection $\mathcal{G}$ and \eqref{eq:Jstability} yield
\begin{align*}
    (\nabla z, \nabla_{\mathrm{pw}}(\mathcal{J} e_h - \mathcal{R} e_h))_{L^2} &= (\nabla_{\mathrm{pw}}(z - \mathcal{G} z), \nabla\mathcal{J} e_h)_{L^2}\\
    &\lesssim \|\nabla_{\mathrm{pw}}(z - \mathcal{G} z)\|_{L^2} \|e_h\|_{a_h}.
\end{align*}
For the second term in \cref{eq:rhsopL21}, using  \cref{eq:defGh} again together with \cref{eq:keyidentity}, we obtain
\begin{align*}
    (\nabla z, \nabla_{\mathrm{pw}} \mathcal{R} e_h)_{L^2} &= (\nabla z, \nabla_{\mathrm{pw}} \mathcal{G} u)_{L^2} - (\nabla_{\mathrm{pw}} \mathcal{G} z, \nabla_{\mathrm{pw}} \mathcal{R} u_h)_{L^2} \\
    &= (\nabla_\pw \mathcal{G}z, \nabla u)_{L^2} - (\nabla_{\mathrm{pw}} \mathcal{G} z, \nabla_{\mathrm{pw}} \mathcal{R} u_h)_{L^2}\\
    & = (\nabla\mathcal{J}\mathcal{I}z,\nabla u)_{L^2} - (\nabla_{\mathrm{pw}} \mathcal{R}\mathcal{I} z, \nabla_{\mathrm{pw}} \mathcal{R} u_h)_{L^2}\\
    & \ \ \ \ \ \ +(\nabla_\pw(\mathcal{G} z - \mathcal{J} \mathcal{I} z), \nabla u)_{L^2}.
\end{align*}
Finally, for the third term in \cref{eq:rhsopL21}, we expand
\begin{align*}
    (V z, \mathcal{J} e_h)_{L^2} = (Vz,\mathcal{J}\mathcal{I}u -\mathcal{J}u_h)_{L^2}.
\end{align*}

To reformulate the sum of the second and third terms in \cref{eq:rhsopL21}, we test the weak formulation \cref{eq:source_pro} with $\mathcal{J}\mathcal{I} z$ and the discrete formulation \cref{eq:dis_source_pro} with $\mathcal{I} z$. Using the resulting identities to rewrite $(\nabla\mathcal{J}\mathcal{I} z,\nabla u)_{L^2}$ and $(\nabla_{\mathrm{pw}}\mathcal{R}\mathcal{I} z,\nabla_{\mathrm{pw}}\mathcal{R}u_h)_{L^2}$~gives
\begin{align}
    &(\nabla z, \nabla_{\mathrm{pw}} \mathcal{R} e_h)_{L^2} + (V z, \mathcal{J} e_h)_{L^2} \notag\\
    \begin{split}\label{eq:rhsproofoptL2}
    &= (f, \mathcal{J} \mathcal{I} z - \Pi_{\mathcal{T}}^{k+1} z)_{L^2} + s_h(\mathcal{I} z, u_h) + (\nabla_{\mathrm{pw}}(\mathcal{G} z - \mathcal{J} \mathcal{I} z), \nabla u)_{L^2}\\
    &\quad + (V u_{\mathcal{T}}, \Pi_{\mathcal{T}}^{k+1} z)_{L^2} - (V \mathcal{J} u_h, z)_{L^2} + (V\mathcal{J}\mathcal{I}u,z)_{L^2} - (Vu,\mathcal{J}\mathcal{I}z)_{L^2}\\
    &\quad - \alpha \left( h_{\mathcal{T}} Q(V u_{\mathcal{T}}), h_{\mathcal{T}} Q(V \Pi_{\mathcal{T}}^{k+1} z) \right)_{L^2}.
    \end{split}
\end{align}

For the first term on the right-hand side of the latter identity, applying the $L^2$-orthogonality of $Q$, the moment-preserving properties of $\mathcal{J}$ (\cref{lem:enrichment}), the approximation estimate~\cref{eq:Japprox} and \cref{eqn:CP} yields
\begin{align}
\label{eq:rhsdefect}
    (f, \mathcal{J} \mathcal{I} z - \Pi_{\mathcal{T}}^{k+1} z)_{L^2}
    &= (Q f, \mathcal{J} \mathcal{I} z - z)_{L^2} + (Q f, Q z)_{L^2} \notag\\
    &\lesssim \mathrm{osc}_{k+1}(f, \mathcal{T}) \|\nabla_{\mathrm{pw}}(z - \mathcal{G} z)\|_{L^2}.
\end{align}

The second term in \cref{eq:rhsproofoptL2} is estimated using the global version of \cref{eq:proof_stabilization}, yielding
\begin{equation*}
    s_h(\mathcal{I} z, u_h) \le |\mathcal{I} z|_{s_h} |u_h|_{s_h} \lesssim \|\nabla_{\mathrm{pw}}(z - \mathcal{G} z)\|_{L^2} |u_h|_{s_h}.
\end{equation*}

To estimate the third term in \cref{eq:rhsproofoptL2}, we note that elementwise integration by parts, together with the moment-preserving properties of $\mathcal{J}$ (cf.~\cref{lem:enrichment}), gives $(\nabla_{\mathrm{pw}} (z - \mathcal{J} \mathcal{I} z), \nabla_{\mathrm{pw}} \mathcal{G} u)_{L^2} = 0$. Combining this with \cref{eq:defGh}, the triangle inequality and the approximation property \cref{eq:Japprox}, we obtain
\begin{align*}
    (\nabla_{\mathrm{pw}}(\mathcal{G} z - \mathcal{J} \mathcal{I} z), \nabla u)_{L^2}
    &= (\nabla_{\mathrm{pw}}(\mathcal{G} z - \mathcal{J} \mathcal{I} z), \nabla_{\mathrm{pw}}(u - \mathcal{G} u))_{L^2} \\
    &\lesssim \|\nabla_{\mathrm{pw}}(z - \mathcal{G} z)\|_{L^2} \|\nabla_{\mathrm{pw}}(u - \mathcal{G} u)\|_{L^2}.
\end{align*}

For the next two terms in \cref{eq:rhsproofoptL2}, the $L^2$-orthogonality properties of $Q$ yield
\begin{equation*}
    (V u_{\mathcal{T}}, \Pi_{\mathcal{T}}^{k+1} z)_{L^2} - (V \mathcal{J} u_h, z)_{L^2}
    = -(Q V u_{\mathcal{T}}, Q z)_{L^2} - (\mathcal{J} u_h - u_{\mathcal{T}}, V z)_{L^2},
\end{equation*}
where, for the first term on the right-hand side, we have
\begin{align*}
    -(Q V u_{\mathcal{T}}, Q z)_{L^2}
    &= -(Q V u, Q z)_{L^2} + (Q V Q u, Q z)_{L^2} + (Q V e_{\mathcal{T}}, Q z)_{L^2} \\
    &\lesssim \Big( \mathrm{osc}_{k+1}(V u, \mathcal{T}) + h_{\max}^{2 - n/p} \|\nabla_{\mathrm{pw}}(u - \mathcal{G} u)\|_{L^2} \\
    &\quad\quad + h_{\max}^{\min\{1, 2 - n/p\}} \|e_h\|_{a_h} \Big) \|\nabla_{\mathrm{pw}}(z - \mathcal{G} z)\|_{L^2},
\end{align*}
which follows from \cref{lem:local_pro_err,lem:higher-order-defect}, and \cref{eqn:CP}.
For the second term on the right-hand side, we obtain
\begin{align*}
    (\mathcal{J} u_h - u_{\mathcal{T}}, V z)_{L^2}
    &= -(\mathcal{J} e_h - e_{\mathcal{T}}, V z)_{L^2} + (\mathcal{J} \mathcal{I} u - \Pi_{\mathcal{T}}^{k+1} u, V z)_{L^2} \\
    &\lesssim \left( h_{\max} \| e_h\|_{a_h}+ h_{\max} \|\nabla_{\mathrm{pw}}(u - \mathcal{G} u)\|_{L^2} \right) \|\nabla z\|_{L^2},
\end{align*}
where we have used the identity $\mathcal{J} e_h - e_{\mathcal{T}} = Q (\mathcal{J} e_h)$, and the stability estimate \cref{eq:Jstability}, as well as the multiplier estimate \cref{eq:continuous-multiplier}, the approximation estimate \cref{eq:Japprox}, and estimate~\cref{eqn:CP}.

The next two terms in \cref{eq:rhsproofoptL2} can be decomposed as
\begin{equation*}
    (V\mathcal{J}\mathcal{I}u,z)_{L^2} - (Vu,\mathcal{J}\mathcal{I}z)_{L^2} = (\mathcal{J}\mathcal{I}u-u,Vz)_{L^2} + (Vu,z-\mathcal{J}\mathcal{I}z)_{L^2}.
\end{equation*}
For the first term, applying the multiplier estimate \cref{eq:continuous-multiplier} together with \cref{eqn:CP} and \cref{eq:Japprox} gives
\begin{equation*}
    (\mathcal{J}\mathcal{I}u-u,Vz)_{L^2} \lesssim h_{\max}\|\nabla_{\mathrm{pw}}(u-\mathcal{G}u)\|_{L^2}\|\nabla z\|_{L^2},
\end{equation*}
while, for the second term, exploiting the moment-preserving properties of $\mathcal{J}$ (cf.~\cref{lem:enrichment}) alongside \cref{eqn:CP} and \cref{eq:Japprox} yields
\begin{align*}
    (Vu,z-\mathcal{J}\mathcal{I}z)_{L^2} &= (Q(Vu),z-\mathcal{J}\mathcal{I}z)_{L^2}
    \lesssim \mathrm{osc}_{k+1}(V u, \mathcal{T})\|\nabla_{\mathrm{pw}} (z-\mathcal{G}z)\|_{L^2}.
\end{align*}

Before estimating the final term in \cref{eq:rhsproofoptL2}, we establish the auxiliary bound
\begin{align*}
    \|h_{\mathcal{T}} Q(V \Pi_{\mathcal{T}}^{k+1} z)\|_{L^2}
    &\le \|h_{\mathcal{T}} Q(V Q z)\|_{L^2} + \|h_{\mathcal{T}} Q(V z)\|_{L^2}\\
    &\lesssim h_{\max}^{2-n/p}\|\nabla_{\mathrm{pw}} (z-\mathcal{G} z)\|_{L^2} + h_{\max}\|\nabla z\|_{L^2},
\end{align*}
which follows from \cref{lem:local_pro_err},  estimate \cref{eq:continuous-multiplier}, and the $L^2$-stability of $Q$.
Combining this bound with \cref{lem:higher-order-defect,lem:local_pro_err}, we obtain for the final term in \cref{eq:rhsproofoptL2}:
\begin{align*}
    &\left( h_{\mathcal{T}} Q(V u_{\mathcal{T}}), h_{\mathcal{T}} Q(V \Pi_{\mathcal{T}}^{k+1} z) \right)_{L^2} \\
    &\quad = \left( h_{\mathcal{T}} Q(V u), h_{\mathcal{T}} Q(V \Pi_{\mathcal{T}}^{k+1} z) \right)_{L^2}
    - \left( h_{\mathcal{T}} Q(V Q u), h_{\mathcal{T}} Q(V \Pi_{\mathcal{T}}^{k+1} z) \right)_{L^2} \\
    &\qquad - \left( h_{\mathcal{T}} Q(V e_{\mathcal{T}}), h_{\mathcal{T}} Q(V \Pi_{\mathcal{T}}^{k+1} z) \right)_{L^2} \\
    &\quad \lesssim \Big( \mathrm{osc}_{k+1}(V u, \mathcal{T}) + h_{\max}^{2 - n/p} \|\nabla_{\mathrm{pw}}(u - \mathcal{G} u)\|_{L^2} + h_{\max}^{\min\{1, 2 - n/p\}} \|e_h\|_{a_h} \Big) \\
    &\qquad \times \Big( h_{\max} \|\nabla z\|_{L^2} + h_{\max}^{2 - n/p} \|\nabla_{\mathrm{pw}}(z - \mathcal{G} z)\|_{L^2} \Big).
\end{align*}
    
Finally, using standard finite element approximation properties together with the dual elliptic regularity estimate~\cref{eq:ellipticregularity}, we obtain
\begin{equation}
\label{eq:estG}
    \|\nabla_{\mathrm{pw}}(z-\mathcal{G} z)\|_{L^2}
    \lesssim h_{\max}^s \|z\|_{H^{1+s}}
    \lesssim h_{\max}^s \|e_{\mathcal{T}}\|_{L^2}.
\end{equation}
Inserting the individual estimates for the terms in~\cref{eq:rhsproofoptL2} into~\cref{eq:rhsopL21}, and combining the resulting bound with \cref{lem:err_est:source,eq:estG} and $\|\nabla z\|_{L^2}\lesssim\|e_{\mathcal T}\|_{L^2}$, yields the assertion.
\end{proof}

\end{document}